\documentclass[proc]{edpsmath}

\pdfoutput = 1

\usepackage[latin1]{inputenc}
\usepackage{array}
\usepackage{amsmath,amsfonts,amssymb,amsthm}
\usepackage{caption}
\usepackage{cite}
\usepackage{color}
\usepackage{dsfont}
\usepackage{epsfig}
\usepackage{float}
\usepackage{graphicx}
\usepackage{enumerate}

\usepackage{amsmath}

 \newcommand{\eg}{\textit{e.g.}}

\newcommand{\Frac}[2]{\displaystyle{\frac{#1}{#2}}}
\newcommand{\abs}[1]{\left|{#1}\right|}
\newcommand{\norm}[1]{\left\|#1\right\|}
\newcommand{\scal}[2]{\left\langle{#1},{#2}\right\rangle}
\newcommand{\tr}{\ensuremath{\mathsf T}}

\newcommand{\OO}{\operatorname{O}}

\newcommand{\set}[2]{\left\{{#1}:{#2}\right\}}
\newcommand{\F}{\ensuremath{\mathcal{F}}}
\newcommand{\X}{\ensuremath{\boldsymbol{\mathcal{X}}}}

\newcommand{\f}{\ensuremath{\textup f}}
\newcommand{\w}{\ensuremath{\textup w}}
\newcommand{\ff}{\ensuremath{\mathbf{f}}}
\newcommand{\ppsi}{\ensuremath{\boldsymbol{\psi}}}
\newcommand{\ww}{\ensuremath{\mathbf{w}}}
\newcommand{\xx}{\ensuremath{\mathbf{x}}}
\newcommand{\yy}{\ensuremath{\mathbf{y}}}
\newcommand{\0}{\ensuremath{\mathbf{0}}}
\newcommand{\1}{\ensuremath{\mathbf{1}}}

\newcommand{\FF}{\ensuremath{\mathbf{F}}}
\newcommand{\PP}{\ensuremath{\mathbf{P}}}
\newcommand{\PPsi}{\ensuremath{\mathbf{\Psi}}}
\newcommand{\WW}{\ensuremath{\mathbf{W}}}
\newcommand{\XX}{\ensuremath{\mathbf{X}}}

\newcommand{\minimize}[3]{\ensuremath{
\begin{aligned}
	\underset{\substack{#2}}{\text{Minimize}} & \quad {#1}
	\\\text{subject to} & \quad {#3}
\end{aligned}
}}

\DeclareMathOperator*{\argmin}{arg\,min}
\DeclareMathOperator*{\argsup}{arg\,sup}

\definecolor{labelkey}{rgb}{0,0.08,0.45}
\definecolor{refkey}{rgb}{0,0.6,0.0}
\definecolor{Brown}{rgb}{0.45,0.0,0.05}
\definecolor{dgreen}{rgb}{0.00,0.49,0.00}
\definecolor{dblue}{rgb}{0,0.08,0.75}
\RequirePackage[colorlinks,hyperindex]{hyperref}
\hypersetup{linktocpage=true,citecolor=dblue,linkcolor=dgreen}

\begin{document}

\title{Emulators for stochastic simulation codes}
\author{Vincent Moutoussamy}
\address{EDF R\&D, MRI, Chatou, France; \email{vincent.moutoussamy@edf.fr}}
\secondaddress{Universit\'e Paul Sabatier, Toulouse, France}
\author{Simon Nanty}
\address{CEA, DEN, DER, F-13108 Saint-Paul-Lez-Durance, France; \email{simon.nanty@cea.fr}}
\secondaddress{Universit\'e Joseph Fourier, Grenoble, France}
\author{Beno\^it Pauwels}
\sameaddress{2}
\secondaddress{IFPEN, DTIMA, Rueil-Malmaison, France; \email{benoit.pauwels@ifpen.fr}}

\begin{abstract} 
Numerical simulation codes are very common tools to study complex phenomena, but they are often time-consuming and considered as \emph{black boxes}. For some statistical studies (\eg{} asset management, sensitivity analysis) or optimization problems (\eg{} tuning of a molecular model), a high number of runs of such codes is needed. Therefore it is more convenient to build a fast-running approximation -- or metamodel -- of this code based on a design of experiments. The topic of this paper is the definition of metamodels for \emph{stochastic} codes. Contrary to deterministic codes, stochastic codes can give different results when they are called several times with the same input. In this paper, two approaches are proposed to build a metamodel of the probability density function of a stochastic code output. The first one is based on \emph{kernel regression} and the second one consists in decomposing the output density on a basis of well-chosen probability density functions, with a metamodel linking the coefficients and the input parameters. For the second approach, two types of decomposition are proposed, but no metamodel has been designed for the coefficients yet. This is a topic of future research. These methods are applied to two analytical models and three industrial cases.
\end{abstract}

\begin{resume} 
Les codes de simulation num\'erique sont couramment utilis\'es pour \'etudier des ph\'enom\`enes complexes. Ils sont cependant souvent co\^uteux en temps de calcul. Pour certaines \'etudes statistiques (\eg{} gestion d'actifs, analyse de sensibilit\'e) ou probl\`emes d'optimisation (\eg{} r\'eglage des param\`etres d'un mod\`ele mol\'eculaire), un grand nombre d'appels \`a de tels codes est n\'ecessaire. Il est alors plus appropri\'e d'\'elaborer une approximation -- ou \emph{m\'etamod\`ele} -- peu co\^uteuse en temps de calcul \`a partir d'un plan d'exp\'erience. Le sujet de cet article est la d\'efinition de m\'etamod\`eles pour des codes \emph{stochastiques}. Contrairement aux codes \emph{d\'eterministes}, les codes stochastiques peuvent renvoyer des r\'esultats diff\'erents lorsqu'ils sont appel\'es plusieurs fois avec le m\^eme jeu de param\`etres. Dans cet article deux approches sont propos\'ees pour la construction d'un m\'etamod\`ele de la densit\'e de probabilit\'e de la sortie d'un code stochastique. La premi\`ere repose sur la \emph{r\'egression \`a noyau} et la seconde consiste \`a d\'ecomposer la densit\'e de la sortie du code dans une base de densit\'es de probabilit\'e bien choisies, avec un m\'etamod\`ele exprimant les coefficients en fonction des param\`etres d'entr\'ee. Pour la seconde approche, deux types de d\'ecomposition sont propos\'es, mais aucun m\'etamod\`ele n'a encore \'et\'e \'elabor\'e pour les coefficients: c'est un sujet de recherches futures. Ces m\'ethodes sont appliqu\'ees \`a deux cas analytiques et trois cas industriels.
 \end{resume}

\maketitle

\section*{Introduction}
	\label{section:intro}
	Numerical simulation codes are very common tools to study complex phenomena. However these computer codes can be intricate and time-consuming: a single run may take from a few hours to a few days. Therefore, when carrying out statistical studies (\textit{e.g.} sensitivity analysis, reliability methods) or optimization processes (\textit{e.g.} tuning of a molecular model), it is more convenient to build an approximation -- or \emph{metamodel} -- of the computer code. Since this metamodel is supposed to take less time to run than the actual computer code, it is substituted for the numerical code. A numerical code is said \emph{deterministic} if it gives the same output each time it is called with a given set of input parameters. For this case, the design of metamodels has already been widely studied.
The present work deals with the formulation of a metamodel for a numerical code which is not deterministic, but \emph{stochastic}.
Two types of numerical codes can be covered by this framework. 
First, the code under consideration may be \emph{intrinsically} stochastic: if the code is called with the same input several times, the output may differ. Thus, the output conditionally to the input is itself a random variable. Second, the computer code may be deterministic with complex inputs. For example, these inputs may be functions or spatiotemporal fields.
These complex inputs are not treated explicitly and are called \emph{uncontrollable parameters}, while the other parameters are considered as \emph{controllable parameters}. When the code is run several times with the same set of controllable parameters, the uncontrollable parameters may change, and thus different output values may be obtained: the output is a random variable conditionally to the controllable parameters.


Two types of approaches can be found in the literature for stochastic code metamodeling.
In the first one, strong hypotheses are made on the shape of the conditional distribution of the output.
In \cite{Rigby2000} the authors consider that the output is normally distributed and estimate its mean and variance, in \cite{Reich2012} the authors assume that the distribution of the output of the code is a mixture of normal distributions and estimate its mean.
In the second approach, the mean and variance of the output are fitted jointly.
In \cite{Zabalza2004,Iooss2009} the authors use joint generalized linear models and joint generalized additive models respectively to fit these quantities.
Joint metamodels based on Gaussian processes have been developed in \cite{Ankenman2010,Marrel2012}.
As we can observe, so far in the literature the metamodels proposed only provide an estimate for one or two of the output distribution moments, sometimes under strong hypotheses on the distribution of the output.

Our approach does not require such a priori hypotheses.
Furthermore we study the whole distribution of the output, allowing to derive estimations for the mean, variance and other moments, and quantiles as well.
The metamodeling is not carried out on the stochastic computer code under consideration itself, but rather on the code $G$ that has same input and that returns the probability density function of the computer code output rather than just a single random value:
\begin{align*}
	G:\X\subset\xR^d & \to\F
	\\\xx & \mapsto G(\xx,\bullet),
\end{align*}
where $\F$ is the space of probability density functions with support in an interval $I$ of $\xR$:
\begin{equation*}
	\F=\set{f:I\to\xR_+}{\int\limits_If=1}.
\end{equation*}
Let $\XX_N=\set{(\xx_i,f_i)}{i=1,\dots,N}\subset(\X\times\F)^N$ be a \emph{training set} ($N\in\xN\setminus\{0\}$): we suppose that we have access to the probability density functions $f_i$ of the outputs associated with $N$ given sets of input parameters $\xx_i$, \ie{}
\begin{equation*}
	f_i=G(\xx_i,\bullet)\qquad(i=1,\ldots, N).
\end{equation*}
In practice, for a given $i$ in $\{1,\ldots, N\}$, the output of a stochastic computer code is not a probability density function $f_i$. We only have access to a finite sample of output values $y_{i1},\ldots, y_{iM}$ which can be written as a vector $\begin{bmatrix}y_{i1} & \cdots & y_{iM}\end{bmatrix}$ representing the function.
In real applications (\cf{} Sections \ref{section:cea}, \ref{section:edf} and \ref{section:ifp}) we chose to derive each probability density function $f_i$ thanks to kernel density estimation.
In this non-parametric method developed by Rosenblatt \cite{Rosenblatt1956} and Parzen \cite{Parzen1962}, the estimator of the probability density function $f_i$, $i=1,\dots,N$ at a point $y$ in $I$ is defined as follows:
\begin{equation*} 
	\hat{f}_i(y)=\frac{1}{M}\sum_{j=1}^MK_H\left(y_{ij}-y\right),
\end{equation*}
where $H\in\mathbb{R}$ is the bandwidth and $K_H$ is the kernel.
Let $\xx_0\in\X$ be a new set of input parameters, we want to estimate the probability density function $f_0$ of the corresponding output:
\begin{align*}
	 f_0=G(\xx_0,\bullet):I\subset\xR & \to\xR_+
	\\t & \mapsto f_0(t)=G(\xx_0,t).
\end{align*}

\begin{figure}
	\center
	\includegraphics[scale=0.4, bb=0 0 956 512]{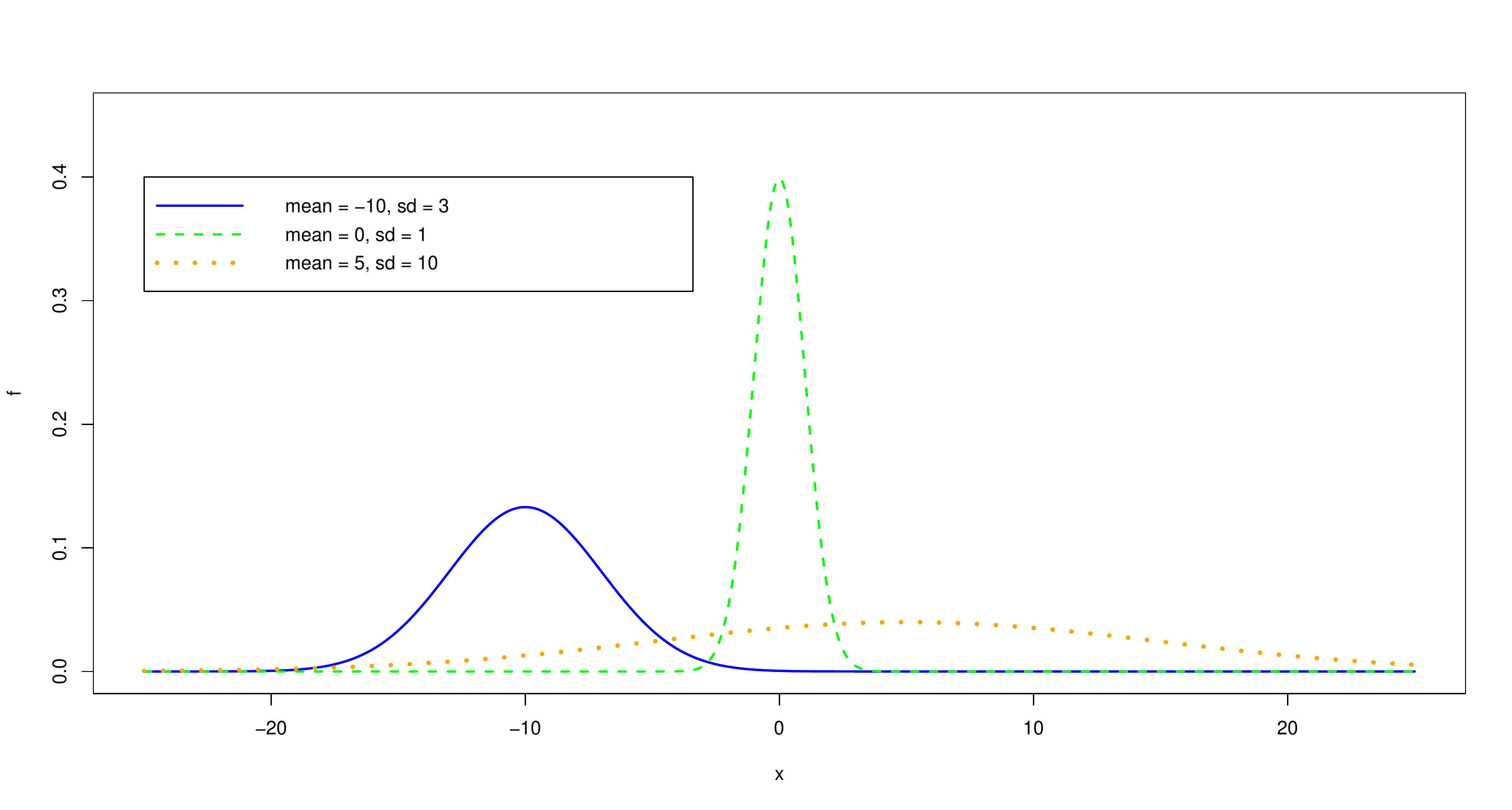}
	\caption{Probability density function given by the numerical code $G$ for different values of the input parameters.}
	\label{figure:gauss:1}
\end{figure}

For example, let us suppose that $\X=\xR\times\xR_+^\ast$, $\xx_0=(\mu_0,\sigma_0)$ and $ G((\mu_0,\sigma_0),\bullet)=\Frac{e^{-\frac{1}{2}\left(\frac{\bullet-\mu_0}{\sigma_0}\right)^2}}{\sqrt{2\pi\sigma_0^2}}. $
Figure \ref{figure:gauss:1} shows the output of $G$ for different values of $\mu_0$ and $\sigma_0$.

Section \ref{Chapter:theory} expounds two types of metamodeling for stochastic computer codes: the first one provides a metamodel with functional outputs based on \emph{kernel regression}, the second one consists in choosing relevant probability density functions to express the probability density of the output as a convex combination of them, the coefficients being linked to the input parameters by another metamodel.
In Section \ref{chapter:numerical_tests}, the previous methods are applied to two analytical test cases.
Finally, the same methods are carried out on three industrial applications in Section \ref{chapter:industrial_applications}.


\section{Metamodels for probability density functions}
	\label{Chapter:theory}
	In this section we expound our two approaches of formulating a metamodel for $G$.
	The first relies on \emph{kernel regression} (\cf{} Section \ref{section:kernel_regression}).
	The second consists in writing $f_0$ as a convex combination of well-chosen probability density functions (\cf{} Section \ref{section:DR}).
	In the following we denote $\norm{\bullet}_{\xLn{2}}$ and $\scal{\bullet}{\bullet}_{\xLn{2}}$ the norm and inner product of $\xLn{2}(I)$ respectively: for all $u$ and $v$ in $\xLn{2}(I)$,
	\begin{align*}
		\norm{u}_{\xLn{2}} = \left(\int\limits_I u^2(t)\xdif{t}\right)^{1/2}\enskip\mathrm{and}\enskip\scal{u}{v}_{\xLn{2}} = \int\limits_I u(t)v(t)\xdif{t}.
	\end{align*}

	\subsection{Kernel regression-based metamodel}
		\label{section:kernel_regression}
		In this subsection, first we apply the classical kernel regression method to the problem under consideration, then a new kernel regression estimator involving the Hellinger distance is introduced and, finally, some perspectives of improvement are suggested.

\subsubsection{Classical kernel regression}
	\label{section:KernelRegression}
\begin{figure}
		\center
		\includegraphics[width = \textwidth, bb=0 0 956 512]{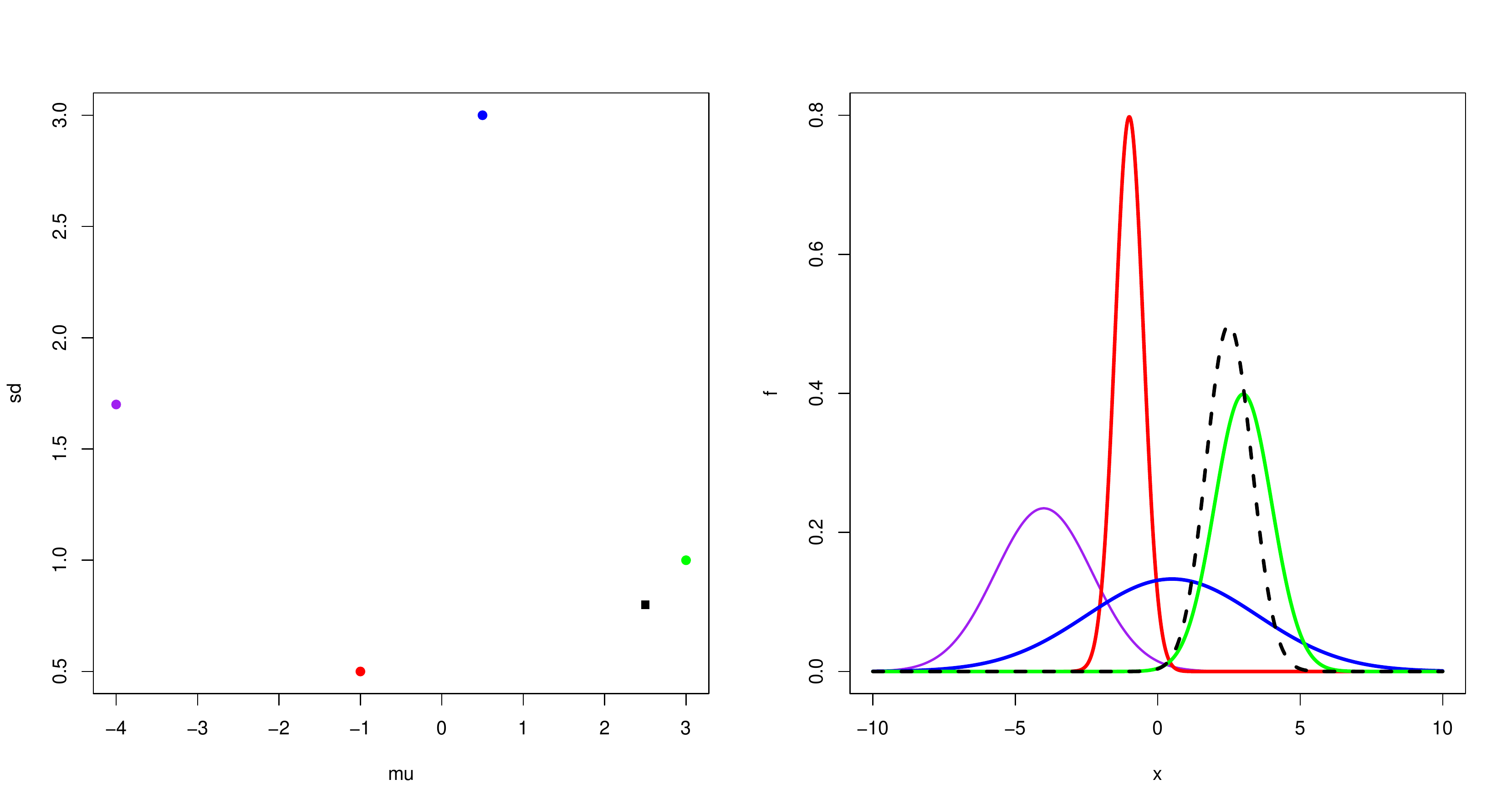}
		\caption{Left: the spots correspond to design points in $\X$. Right: the associated probability density functions given by the code $G$ defined in the introduction. The probability density function associated to the black square (dashed line) in left Figure seems to be closer to the probability density function represented in green than to the others.}
		\label{figure:gauss:2}
	\end{figure}

	Let $\xx_0\in\X$ and $f_0=G(\xx_0,\bullet)$. We suppose that the sample set $\XX_N=\set{\left(\xx_i,f_i\right)}{i=1,\dots,N}\subset(\X\times\F)^N$ is available. We want to estimate $f_0$ by $\hat{f}_0\in\F$, knowing $\XX_N$. However, forcing $\hat{f}_0$ to be in $\F$ can be difficult in practice. We propose a first estimator given by
	\begin{equation*}
		\hat{f}_0=\sum\limits_{i=1}^N\alpha_if_i
	\end{equation*}
	where $\alpha_i\in\xR$ ($i=1,\ldots, N$). In order to have $\hat{f}_0$ in $\F$, we impose
	\begin{align*}
		\alpha_i & \geq0\quad\left(i=1,\ldots, N\right),\\
		\sum\limits_{i=1}^N\alpha_i & =1.
	\end{align*}
	The \emph{non-parametric kernel regression} verifies the previous constraints. This estimator, introduced in \cite{nadaraya1964,watson1964}, can be written as follows:
	\begin{equation}
		\label{eq:KR}
 		\hat{f}_0=\sum\limits_{i=1}^N\Frac{K_H(\xx_i,\xx_0)}{\sum\limits_{j=1}^NK_H(\xx_j,\xx_0)}f_i,
	\end{equation}
	where $K_H:\xR^d \times \xR^d \rightarrow \xR$ is a \emph{kernel function}. In the following, the \emph{Gaussian kernel} is used:
	\begin{equation*}
 		K_H(\xx,\yy)= \frac{1}{\sqrt{2\pi\det(H)}}e^{ -(\xx-\yy)^\tr H^{-1}(\xx-\yy)}\quad\left(\xx,\yy\in\xR^d\right),
	\end{equation*}
	where
	\begin{eqnarray*}
 		H=\operatorname{diag}(h_1,\ldots, h_d) 
	\end{eqnarray*}
	is the (diagonal) \emph{bandwidth matrix}, with $h_1,\ldots, h_d\in\xR_+$. The higher $h_j$ is, the more points $\xx_i$ are taken into account in the regression in the direction $j$. The method relies on the intuition that if $\xx_0$ is close to $\xx_i$ then $f_0$ will be more influenced by the probability density function $f_i$ (see Figure \ref{figure:gauss:2}).
	We want that the bandwidth matrix $H^\ast$ minimizes the global error of estimation. That means finding $H^\ast=\operatorname{diag}(h_1^\ast,\dots,h_d^\ast)$ such that
	\begin{equation*}
		\operatorname{diag}(h_1^\ast,\ldots,h_d^\ast)\in\argmin_{h_1,\ldots,h_d\in\xR_+}\int_{\X}\norm{\hat{f}_0-f_0}_{\xLn{2}}^2\xdif{\xx_0}.
	\end{equation*}
	In order to get an approximation of $H^*$, one uses the \emph{leave-one-out cross-validation}, as proposed in Hardle and Marron \cite{Hardle1985}. Then
	\begin{equation*}
 		H^{\ast}\in\argmin_{h_1,\ldots,h_d\in\xR_+}\sum\limits_{i=1}^N\norm{\hat{f}_{-i}-f_i}_{\xLn{2}}^2,
	\end{equation*}
	with
	\begin{align*}
		\hat{f}_{-i} & =\sum\limits_{\substack{j=1\\j\neq i}}^N\alpha_j^{-i}f_j
		\\\alpha_j^{-i}&=\Frac{e^{-(\xx_i-\xx_j)^t H^{-1}(\xx_i-\xx_j)}}{\displaystyle{\sum\limits_{\substack{l=1\\l\neq i}}^Ne^{-(\xx_i-\xx_l)^t H^{-1}(\xx_i-\xx_l)}}}.
	\end{align*}
	This optimization problem is solved with the L-BFGS-B algorithm developed by Byrd et al. \cite{Byrd1995}. Finally, the kernel regression metamodel is given by
	\begin{equation}
		\label{eq:L2:1}
		\hat{f}_0=\sum\limits_{i=1}^N\alpha_{i,H^{\ast}}f_i,
	\end{equation}
	with
	\begin{equation*}
  		\alpha_{i,H^{\ast}} = \frac{ e^{-(\xx_0-\xx_i)^T (H^\ast)^{-1}(\xx_0-\xx_i)}}{\displaystyle{\sum\limits_{j=1}^N e^{-(\xx_0-\xx_j)^T (H^\ast)^{-1}(\xx_0-\xx_j)}}}.
	\end{equation*}

%
%
%

\subsubsection{Kernel regression with the Hellinger distance}
	\label{section:KernelRegressionHellinger}
	The kernel regression introduced in equation \eqref{eq:KR} can be considered as a minimization problem. The estimator of the function $f_0$ corresponding to the point $\xx_0$ is as follows:
	\begin{equation}
		\label{eq:KR2}
		 \hat{f}_0=\argmin_{f\in\F}\sum_{i=1}^NK_H(\xx_i,\xx_0)\int\limits_I\left(f_i-f\right)^2.
	\end{equation}
	The kernel regression is therefore a locally weighted constant regression, which is fitted by minimizing weighted least squares.
	The equivalence between estimators \eqref{eq:KR} and \eqref{eq:KR2} is proven in Appendix \ref{section:L2.kernel.equivalence}.
	We observe that the $\xLn{2}$ distances between the sample functions and the unknown function appear in the objective function of the minimization problem \eqref{eq:KR2}.
	

	Another classical example of distance between probability density functions is the Hellinger distance.
	We derived a new kernel regression estimator with respect to this distance replacing the $\xLn{2}$ distance by the Hellinger distance in equation \eqref{eq:KR2}.
	The kernel estimator thus takes the following form:


	\begin{align}
  		\hat f_0 = \argmin_{f\in\mathcal F} \sum_{i=1}^N K_H(\xx_i,\xx_0)\int\limits_I\left(\sqrt{f_i(t)} - \sqrt{f(t)}\right)^2\xdif{t}.
  		\label{eq:Hellinger:2}
	\end{align}
	The following analytical expression can be derived for this Hellinger distance-based estimator:
	\begin{equation}
		\label{eq:Hellinger:1}
		\hat{f}_0=\Frac{\displaystyle\left(\sum\limits_{i=1}^NK_H(\xx_i,\xx_0)\sqrt{f_i}\right)^2}{\displaystyle\int\limits_I\left(\sum\limits_{i=1}^N K_H(\xx_i,\xx_0)\sqrt{f_i}\right)^2}.
	\end{equation}
	The calculation that leads to \eqref{eq:Hellinger:1} is provided in Appendix \ref{section:Hellinger.kernel.equivalence}.

	\begin{figure}
		\center
		\includegraphics[scale=0.4, bb=0 0 619 417]{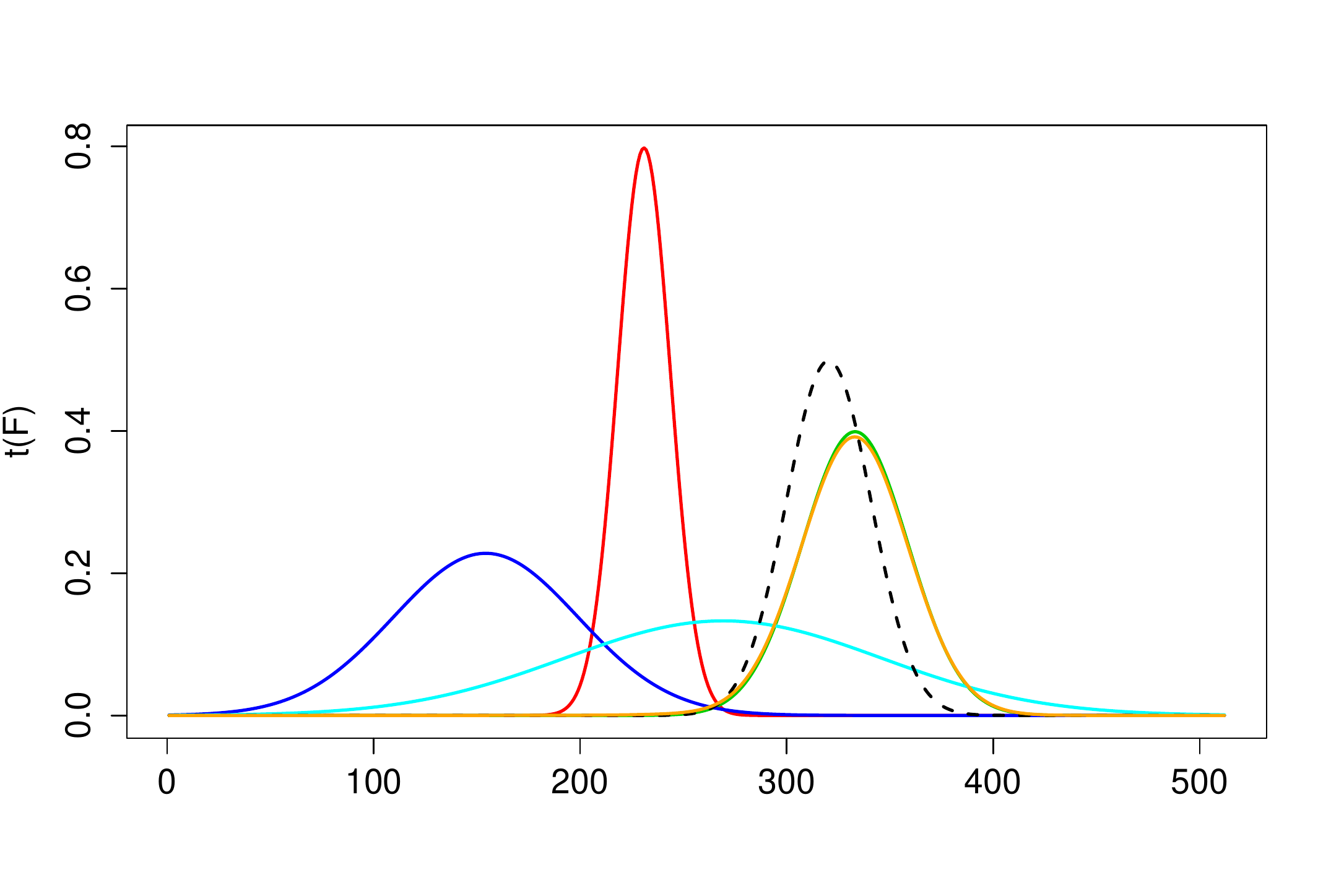}
		\caption{The orange curve is the Hellinger kernel prediction for the (unknown) output in black, built from the dark blue, red, light blue and green learning probability density functions ($N=4$).}
		\label{figure:Hellinger.kernel.example}
	\end{figure}
	
	From now on, we will denote $\hat{f}_{0,\xLn{2}}$ and $\hat{f}_{0,\mathrm{He}}$ the estimators given by equations \eqref{eq:L2:1} and \eqref{eq:Hellinger:1} respectively.
	Let us illustrate the Hellinger kernel regression on the four gaussian density functions given on Figure \ref{figure:gauss:2}.
	On Figure \ref{figure:Hellinger.kernel.example} the same $N=4$ learning curves (in dark blue, red, light blue and green) are represented, along with an unknown output $f_0$ (in dashed black) and its associated Hellinger kernel estimator $\hat{f}_{0,\mathrm{He}}$ (in orange).
	We can observe that, although the output is not very well estimated, most of the weight is given to the curve corresponding to the closest vector of parameters (in green).

\subsubsection{Perspectives}
	
	One of the well-known major drawbacks of the kernel regression is due to the curse of dimensionality.
	The bandwidth estimation quality is deteriorating as the dimension $d$ of the problem increases for a constant sample size.
	First, as the sample points are sparsely distributed in high dimension, the local averaging performs poorly.
	Second, the optimization problem to find the optimal bandwidth is more complex in higher dimension. 	
	In the presented work, kernel regression has been used with a fixed bandwidth for all sample points.
	The Adaptive or "Variable-bandwidth" kernel estimation techniques have been proposed in the literature to use varying bandwidth over the domain of the parameters $\mathcal X$.
	See Muller and Stadtmuller \cite{Muller1987} and Terrell and Scott \cite{Terrell1992} for more details on these methods. Terrell and Scott \cite{Terrell1992} review several adaptive methods and compare their asymptotic mean squared error, while Muller and Stadtmuller \cite{Muller1987} propose a varying bandwidth selection method.
	In numerous cases, using non-constant bandwidths leads to better estimates as it enables to better capture the different features of the curves.
	Hence, the quality of the proposed methods could be improved by using adaptive kernel regression.
	However, the use of kernel regression is not advised in high dimension.
	


	\subsection{Functional decomposition-based metamodel}
		\label{section:DR}
		In this section, we aim at building a set of \emph{basis functions} $w_1,\ldots, w_q$ in $\F$ to approximate $f_0$ by an estimator of the form
\begin{equation}
	\hat{f}_0=\sum\limits_{k=1}^q\psi_{k}(\xx_0)w_k,
	\label{eq:functional.decomposition}
\end{equation}
where $\psi_{k}$ are functions from $\X$ to $\xR$ ($k=1,\ldots, q$), called \emph{coefficient functions}, satisfying
\begin{equation*}
	\left\{
	\begin{aligned}
		&\psi_{k}(\xx)\geq0
		\\&\sum_{k=1}^q\psi_{k}(\xx)=1
	\end{aligned}
	\right.\qquad\left(\xx\in\X\right).
\end{equation*}
Kernel regression and functional decomposition have different goals but have similar forms.
The former provides an estimator which is a convex combination of all the $N$ sampled probability density functions $f_1,\ldots,f_N$.
The latter gives a combination of fewer density functions built from the experimental data, thus lying in a smaller space, in order to \emph{reduce the dimension} of the problem.
We observe that the estimator provided by classical kernel regression \eqref{eq:KR} is of the form \eqref{eq:functional.decomposition} with $q=N$ (no dimension reduction), $w_k=f_k$ for $k=1,\ldots,q$ and
\begin{equation*}
	\psi_k\left(\xx_0\right)=\Frac{K_H(\xx_k,\xx_0)}{\sum\limits_{j=1}^NK_H(\xx_j,\xx_0)}\qquad\left(k=1,\ldots,q\right).
\end{equation*}

We will build the basis functions and the coefficient functions so that they fit the available data $\XX_N=\set{\left(\xx_i,f_i\right)}{i=1,\ldots, N}$. Thus the functions $w_1,\ldots, w_q$ will be designed so that, for all $i$ in $\{1,\ldots, N\}$, $f_i$ is approximated by
\begin{align}
	&\hat{f}_i=\sum\limits_{k=1}^q\psi_{ik}w_k\nonumber
	\\&\mathrm{where}\enskip
	\left\{\begin{aligned}
	&\psi_{ik}=\psi_k(\xx_i)&&(k=1,\ldots, q)
	\\&\psi_{ik}\geq0&&(k=1,\ldots, q)
	\\&\sum_{k=1}^q\psi_{ik}=1.
	\end{aligned}\right.
	\label{eq:psi_ik.constraints}
\end{align}
For all $i$ in $\{1,\ldots, N\}$, we require the approximation $\hat{f}_i$ of $f_i$ to be a probability density function as well; therefore $\hat{f}_i$ must be non-negative and have its integral equal to 1.

We study three ways of choosing the basis functions.
In Section \ref{section:CPCA} we discuss the adaptation of \emph{Functional Principal Components Analysis} (FPCA) to compute the basis functions: $w_1,\ldots,w_q$ are then orthonormal with integral equal to 1 and their negative values are interpreted as 0, hence the predictions on the experimental design are probability density functions.
In the following two subsections we propose two new approaches to probability density function decomposition.
The non-negativity constraint is taken into account directly rather than through a post-treatment of FPCA results (setting negative values to zero and renormalizing).
In Section \ref{section:EIM}, we propose to adapt the \emph{Empirical Interpolation Method} or \emph{Magic Points method} \cite{Maday2007} to our problem: we obtain an algorithm to select the basis functions (not necessarily orthogonal) among the sample distributions $f_1,\ldots,f_N$ (while losing the interpolation property). It follows that the predictions on the experimental design are probability density functions.
In Section \ref{section:AQM} a third approach is proposed, consisting in computing directly the basis functions that minimize the $\xLn{2}$-approximation error on the experimental design (without imposing orthogonality):
\begin{equation}
	\label{eq:L2.approximation.error}
	\sum\limits_{i=1}^N\norm{f_i-\sum\limits_{k=1}^q\psi_{ik}w_k}_{\xLn{2}}^2.
\end{equation}
We introduce an iterative procedure to ease the numerical computation.
The basis functions obtained are probability density functions and so are the predictions on the experimental design.
Finally, in Section \ref{section:coefficients} we discuss the formulation of the coefficient functions.

		\subsubsection{Constrained Principal Component Analysis (CPCA)}
			\label{section:CPCA}
			Among all the \emph{dimension reduction} techniques, one of the best known is the \emph{Functional Principal Component Analysis} (FPCA), developed by Ramsay and Silverman \cite{Ramsay2005}, which is based on \emph{Principal Component Analysis} (PCA).
Let us denote $\bar{f}$ the \emph{mean estimator} of $f_1,\ldots, f_N$:
\begin{equation*}
	\bar{f}=\Frac{1}{N}\sum\limits_{i=1}^Nf_i.
\end{equation*}
The goal of FPCA is to build an orthonormal basis $w_1,\ldots, w_q$ for the \emph{centered functions} $f_i^c=f_i-\bar f$ ($i=1,\ldots, N$) so that the projected variance of every centered function onto the functions $w_1,\ldots, w_q$ is maximized.
The following maximization problem is thus considered:
\begin{equation*}
	\minimize{\sum_{i=1}^N\norm{f^c_{i}-\sum_{k=1}^q\scal{f^c_i}{w_k}_{\xLn{2}}w_k}_{\xLn{2}}^2}{w_1,\ldots,w_q:I\to\xR}{\scal{w_k}{w_l}_{\xLn{2}}=\delta_{kl}\qquad(k,l=1,\ldots,q),}
\end{equation*}
where $\delta_{kl}$ is 1 when $k=l$ and 0 otherwise.
Therefore the estimator of $f_i$ $(i=1,\ldots,N)$ is
\begin{equation*}
	\hat{f}_i=\bar{f}+\sum\limits_{k=1}^q\psi_{ik}w_k,
\end{equation*}
with $\psi_{ik}=\scal{f^c_i}{w_k}_{\xLn{2}}$ ($k=1,\ldots,q$).
To solve this minimization problem, Ramsay and Silverman \cite{Ramsay2005} proposes to express the sample functions on a spline basis.
Then PCA can be applied to the coefficients of the functions on the spline basis.
They also propose to apply PCA directly to the discretized functions.
The FPCA decompostion ensures that $\int_I \hat{f}_i=\int_I f_i$ for $i=1,\dots,N$ so that the integrals of the sample probability density functions approximations $\hat{f}_i$ are equal to one.
However, the FPCA decomposition does not ensure that the estimators are non-negative.

Delicado \cite{Delicado2011} proposes to apply FPCA to the logarithm transform $g_i=\log f_i$ to ensure the positivity of the prediction $\hat{f}_i=\exp\hat{g}_i$ ($i=1,\ldots,N$).
However this approach does not ensure that the approximations are normalized. 

Affleck-Graves et al. \cite{Affleckgraves1979} have proposed to put a non-negativity constraint on the first basis function $w_1$ and let the other basis functions free.
Indeed they claim that in general it is not possible to put such a constraint on the basis functions $w_2, \dots , w_q$.
Hence the first basis function is forced to be non-negative but not the other ones. Therefore, it cannot be guaranteed that the approximations are non-negative.

Another method has been proposed by Kneip and Utikal \cite{Kneip2001} to ensure that the FPCA approximations are non-negative.
The FPCA method is applied to the sample functions.
Then the negative values of the probability density functions approximations are interpreted as 0 and they are normalized to one.
In the numerical studies of Sections \ref{chapter:numerical_tests} and \ref{chapter:industrial_applications}, we will refer to this method as CPCA, for Constrained Principal Component Analysis.

		\subsubsection{Modified Magic Points (MMP)}
			\label{section:EIM}
			The \emph{Empirical Interpolation Method} (EIM) or \emph{Magic points} (MP) \cite{Maday2007} is a \emph{greedy} algorithm that builds interpolators for a set of functions.
Here the set of functions under consideration is $G\left(\X,\bullet\right)=\set{G(\xx,\bullet)}{\xx\in\X}$.
The algorithm iteratively picks a set of basis functions in $G\left(\X,\bullet\right)$ and a set of \emph{interpolation points} in $I$ to build an interpolator with respect to these points.
Let $\xx$ be in $\X$ and $f_\xx=G\left(\xx,\bullet\right)$.
At step $q-1$, the current basis functions and interpolation points are denoted $w_1,\ldots,w_{q-1}:I\to\xR$ and $t_{i_1},\ldots,t_{i_{q-1}}\in I$ respectively.
The interpolator $I_{q-1}[f_\xx]$ of $f_\xx$ is defined as a linear combination of the basis functions:
\begin{equation*}
	I_{q-1}[f_\xx]=\sum\limits_{k=1}^{q-1}\psi_k(\xx)w_k,
\end{equation*}
where the values of the coefficient functions $\psi_1,\ldots,\psi_{q-1}:\X\to\xR$ at $\xx$ are uniquely defined by the following \emph{interpolation equations} (see \cite{Maday2007} for existence and uniqueness):
\begin{equation*}
	I_{q-1}[f_\xx]\left(t_{i_k}\right)=f_{\xx}\left(t_{i_k}\right)\qquad\left(k=1,\ldots,q-1\right).
\end{equation*}
At step $q$, one picks the parameter $\xx_{i_q}$ such that $f_{\xx_{i_q}}\left(=G\left(\xx_{i_q},\bullet\right)\right)$ is the element of $G\left(\X,\bullet\right)$ whose $\xLn{2}$-distance from its own current interpolator $I_{q-1}[f_{\xx_{i_q}}]$ is the highest.
The next interpolation point $t_{i_q}$ is the one that maximizes the gap between $f_{\xx_{i_q}}$ and $I_{q-1}[f_{\xx_{i_q}}]$.
Then the new basis function $w_q$ is defined as a particular linear combination of $f_{\xx_{i_q}}$ and $I_{q-1}[f_{\xx_{i_q}}]$.
Let us be more specific by summarizing the algorithm \cite{Bebendorf2013} hereafter.
\begin{lgrthm}[MP]\
	\begin{itemize}
		\item Set $q \leftarrow 1$ and $I_0\leftarrow 0$.
		\item While $\varepsilon > tol$ do
			\begin{enumerate}[(a)]
				\item choose $\xx_{i_q}$ in $\X$:
					\begin{equation*}
						\xx_{i_q}\leftarrow\underset{\xx\in\X}{\argsup}\norm{f_{\xx}-I_{q-1}[f_{\xx}]}_{\xLn{2}}
					\end{equation*}
					and the associated interpolation point $t_{i_q}$:
					\begin{equation*}
						t_{i_q}\leftarrow\underset{t\in I}{\argsup}\abs{f_{\xx_{i_q}}(t)-I_{q-1}[f_{\xx_{i_q}}](t)},
					\end{equation*}
				\item define the next element of the basis:
					\begin{equation*}
						w_q\leftarrow\Frac{f_{\xx_{i_q}}(\bullet)-I_{q-1}[f_{\xx_{i_q}}](\bullet)}{f_{\xx_{i_q}}(t_{i_q})-I_{q-1}[f_{\xx_{i_q}}](t_{i_q})},
					\end{equation*}
				\item compute the \emph{interpolation error}:
					\begin{equation*}
						\varepsilon\leftarrow\norm{err_q}_{\xLinfty(\X)},\enskip\text{where}\enskip err_q:\xx\mapsto\norm{f_\xx-I_q[f_\xx]}_{\xLn{2}},
					\end{equation*}
				\item set $q\leftarrow q+1$.
			\end{enumerate}
	\end{itemize}
\end{lgrthm}
This method has been successfully applied \eg{} to a heat conduction problem, crack detection and damage assessment of flawed materials, inverse scattering analysis \cite{Cuong2005}, parameter-dependent convection-diffusion problems around rigid bodies \cite{Tonn2006}, in biomechanics to describe a stenosed channel and a bypass configuration \cite{Rozza2009}.

In general, the interpolators provided by the MP algorithm are not probability density functions.
Therefore we modified the MP method so that the approximations are non-negative with integral equal to one.
In the derived new method the interpolation is lost, but the basis functions are picked in a similar greedy way.
In the following we denote $A_{q-1}[f_\xx]$ the interpolator at step $q-1$ of the function $f_\xx$ associated with a parameter $\xx$ in $\X$:
\begin{eqnarray*}
	A_{q-1}[f_\xx]=\sum_{k=1}^{q-1}\psi_k(\xx)w_k,
\end{eqnarray*}
where the values of the coefficient functions at $\xx$ are now defined as solutions of the following convex quadratic program:
\begin{equation}
	\label{eq:minimization.MP}
	\minimize
	{\norm{f_\xx-\sum\limits_{k=1}^{q-1}\psi_k(\xx)w_k}_{\xLn{2}}^2}
	{\psi_1(\xx),\ldots,\psi_{q-1}(\xx)\in\xR}
	{\left\{\begin{aligned}&\psi_k(\xx)\geq0&&(k=1,\ldots,q-1)\\&\sum\limits_{k=1}^{q-1}\psi_k(\xx)=1.&&\end{aligned}\right.}
\end{equation}
At step $q$, the parameter $\xx_{i_q}$ is chosen in the same way as in the MP algorithm. The new basis function $w_q$ is then $f_{\xx_{i_q}}$.
The \emph{Modified Magic Points} (MMP) algorithm is detailed below.
\begin{lgrthm}[MMP]\
	\begin{itemize}
		\item Set $q \leftarrow 1$ and $I_0\leftarrow 0$.
		\item While $\varepsilon > tol$ do
			\begin{enumerate}[(a)]
				\item\label{choice_step} choose $\xx_{i_q}$ in $\X$:
					\begin{equation*}
						\xx_{i_q}\leftarrow\underset{\xx\in\X}{\argsup}\norm{f_{\xx}-A_{q-1}[f_{\xx}]}_{\xLn{2}},
					\end{equation*}
				\item define the next element of the basis:
					\begin{equation*}
						w_q\leftarrow f_{\xx_{i_q}},
					\end{equation*}
				\item compute the \emph{estimation error}:
					\begin{equation*}
						\varepsilon\leftarrow\norm{err_q}_{\xLinfty(\X)},\enskip\text{where}\enskip err_q:\xx\mapsto\norm{f_\xx-A_q[f_\xx]}_{\xLn{2}},
					\end{equation*}
				\item set $q\leftarrow q+1$.
			\end{enumerate}
	\end{itemize}
\end{lgrthm}
\begin{rmrk}
	In step (\ref{choice_step}) of the MMP algorithm, the $\xLn{2}$ norm can be replaced by the Hellinger distance.
	This yields a new algorithm which is tested along with the previous one in Section \ref{chapter:numerical_tests}. 
\end{rmrk}
In practice, we do not apply the MMP algorithm to the whole set $G\left(\X,\bullet\right)$ since the computation of the functions is time-consuming.
It is rather applied to the available sample set $\{f_1,\ldots,f_N\}$.
Thus the MMP algorithm is a way to select the most relevant functions among the sample $f_1,\ldots,f_N$.
Then the step (\ref{choice_step}) of the algorithm is only a finite maximization.
Furthermore we build a regular grid of $M$ points $t_1,\ldots,t_M$ in the interval $I$:
\begin{equation}
	\label{eq:regular.grid}
	\begin{split}
		&t_j=t_1+(j-1)\Delta t\qquad(j=1,\ldots, M)
		\\&\mathrm{with}\enskip\Delta t\in\xR_+,
	\end{split}
\end{equation}
and we discretize the functions $f_i$ ($i=1,\ldots,N$) and $w_k$ ($k=1,\ldots,q$) on this grid:
\begin{equation*}
	\begin{split}f_{ij}&=f_i(t_j)\\\w_{kj}&=w_k(t_j)\end{split}\qquad(j=1,\ldots, M).
\end{equation*}
Thus we can rewrite the problem (\ref{eq:minimization.MP}) in the following way for all $i\in\{1,\ldots,N\}$: 
\begin{equation}
	\label{eq:minimization.MP.discretized}
	\minimize
	{\sum\limits_{j=1}^M\left[f_{ij}-\sum\limits_{k=1}^q\psi_{ik}w_{kj}\right]^2\Delta t}
	{\psi_{i1},\ldots,\psi_{iq}}
	{\left\{\begin{aligned}\psi_{i1},\ldots,\psi_{iq}&\geq0\\\sum\limits_{k=1}^q\psi_{ik}&=1.\end{aligned}\right.}
\end{equation}
These $N$ problems are convex quadratic programs with $q$ unknowns and $q+1$ linear constraints, they are easily solvable for basis sizes smaller than $10^4$.
\begin{rmrk}
	Let us suppose $q=1$, \ie{} the basis contains only one function. The equality constraint in the minimization \eqref{eq:minimization.MP.discretized} reduces to $\psi_{1i}=1$: the approximation is the same for each function $f_i$ in the sample.
\end{rmrk}
%
%

		\subsubsection{Alternate Quadratic Minimizations (AQM)}
			\label{section:AQM}
Our third functional decomposition approach consists in tackling directly the problem consisting in minimizing the $\xLn{2}$ norm of the approximation error \eqref{eq:L2.approximation.error} with basis functions $w_1,\ldots, w_q$ in $\F$ and coefficients $\psi_{ik}$ satisfying the constraints \eqref{eq:psi_ik.constraints}. It differs from PCA because no orthonormality condition is imposed, and the decomposition is carried out directly on the raw data, \ie{} without centering.
This functional minimization program can be written as follows.
\begin{equation}
	\label{functional minimization}
	\minimize
	{\Frac{1}{2}\sum\limits_{i=1}^N\norm{f_i-\sum\limits_{k=1}^q\psi_{ik}w_k}_{\xLn{2}}^2}
	{w_k,\psi_{ik}\\k=1,\ldots, q,i=1,\ldots, N}
	{\left\{\begin{aligned}&w_k\in\F & &(k=1,\ldots, q) &\\&\psi_{ik}\geq0 & &(i=1,\ldots, N)(k=1,\ldots, q)\\&\sum\limits_{k=1}^q\psi_{ik}=1 & &(i=1,\ldots, N).\end{aligned}\right.}
\end{equation}


Just as before, instead of working on actual functions, we consider discretizations of them. We discretize the interval $I$ in the same way \eqref{eq:regular.grid} and we set
\begin{align*}
	\ff_i&=\begin{bmatrix}f_i(t_1) & \cdots & f_i(t_M)\end{bmatrix}\qquad\left(i=1,\ldots, N\right)
	\\&=\begin{bmatrix}\f_{i1} & \cdots & \f_{iM}\end{bmatrix}
	\\\ww_k&=\begin{bmatrix}w_k(t_1) & \cdots & w_k(t_M)\end{bmatrix}\qquad\left(k=1,\ldots, q\right)
	\\&=\begin{bmatrix}\w_{k1} & \cdots & \w_{kM}\end{bmatrix}.
\end{align*}
The functional minimization program (\ref{functional minimization}) is then replaced by the following vectorial minimization program (in this section $\norm{\bullet}_{\xR^M}$ designates the euclidean norm on $\xR^M$).
\begin{equation}
	\label{vectorial minimization}
	\minimize
	{\Frac{1}{2}\sum\limits_{i=1}^N\norm{\ff_i-\sum\limits_{k=1}^q\psi_{ik}\ww_k}_{\xR^M}^2}
	{\ww_k,\psi_{ik}\\k=1,\ldots, q,i=1,\ldots, N}
	{\left\{\begin{aligned}&\w_{kj}\geq0 & &(k=1,\ldots, q)(j=1,\ldots, M)\\&\sum\limits_{j=1}^M\w_{kj}\Delta t=1 & &(k=1,\ldots, q)\\&\psi_{ik}\geq0 & &(i=1,\ldots, N)(k=1,\ldots, q)\\&\sum\limits_{k=1}^q\psi_{ik}=1 & &(i=1,\ldots, N).\end{aligned}\right.}
\end{equation}
Let us introduce additional notations to reformulate the objective function and the variables in a more compact way.
\begin{align*}
	\PPsi=&\left[\psi_{ik}\right]_{i=1,\ldots, N; k=1,\ldots, q}=\begin{bmatrix}\ppsi_1\\\vdots\\\ppsi_N\end{bmatrix}
	&\WW=&\left[\w_{kj}\right]_{k=1,\ldots, q; \ j=1,\ldots, M}=\begin{bmatrix}\ww_1\\\vdots\\\ww_q\end{bmatrix}
\end{align*}
We recall that the Frobenius norm of a matrix $\mathbf{A}=\left[a_{ij}\right]_{i=1,\ldots, N; j=1,\ldots, M}\in\xR^{N\times M}$ is defined as
\begin{equation*}
	\norm{\mathbf{A}}_{\text{F}}=\sqrt{\sum\limits_{i=1}^N\sum\limits_{j=1}^Ma_{ij}^2}.
\end{equation*}
This being said, the vectorial minimization program \eqref{vectorial minimization} can be written in matrix form.
\begin{equation}
	\label{matricial minimization}
	\minimize{\OO\left(\PPsi,\WW \right)=\Frac{1}{2}\norm{\FF-\PPsi\WW}_{\text{F}}^2}{\PPsi\in\xR_+^{N\times q},\WW\in\xR_+^{q\times M} }
	{\left\{\begin{aligned}& \ppsi_i\1=1 &(i=1,\ldots, N)\\&\ww_k\1=1/\Delta t &(k=1,\ldots, q).\end{aligned}\right.}
\end{equation}

The objective function of the program \eqref{vectorial minimization} is second-order polynomial, hence twice continuously differentiable (but not necessarily convex). Its first-order derivatives are expressed as
\begin{align*}
	\partial_{\psi_{ik}}\OO\left(\PPsi,\WW\right)&=-\ww_k\left(\ff_i-\ppsi_i\WW\right)^\tr
	&&(i=1,\ldots, N)(k=1,\ldots, q),
	\\\partial_{w_{kj}}\OO\left(\PPsi,\WW\right)&=-\PPsi_{\bullet,k}^\tr\left(\FF_{\bullet,j}-\PPsi\WW_{\bullet,j}\right)
	&&\left(k=1,\ldots, q\right)\left(j=1,\ldots, M\right),
\end{align*}
where $\PPsi_{\bullet,k}$ is the $k^{\mathrm{th}}$ column of $\PPsi$, $\FF_{\bullet,j}$ and $\WW_{\bullet,j}$ are the $j^{\mathrm{th}}$ columns of $\FF$ and $\WW$ respectively, and here are the second-order derivatives which will be useful further down:
\begin{align*}
	\partial^2_{\psi_{ik}\psi_{ik'}}\OO\left(\PPsi,\WW\right)&=\ww_k\ww_{k'}^\tr
	&&(i=1,\ldots, N)(k,k'=1,\ldots, q),
	\\\partial^2_{w_{kj}w_{kj'}}\OO\left(\PPsi,\WW\right)&=\delta_{jj'}\norm{\PPsi_{\bullet,k}}_{\xR^N}^2
	&&\left(k=1,\ldots, q\right)\left(j,j'=1,\ldots, M\right).
\end{align*}

The program \eqref{matricial minimization} has $q(M+N)\approx10^4$ design variables, $q(M+N)+q+N\approx10^4$ constraints and may not be convex. Hence its resolution through the use of a numerical solver may be very time-consuming (it turns out we could not manage to get a satisfactory feasible solution). In order to circumvent this issue we chose to implement \eqref{matricial minimization} as successive convex quadratic minimization programs. The value of the program can be written as follows.
\begin{equation*}
	\inf\limits_{\substack{\left(\PPsi,\WW\right)\in\left(\xR_+^{N\times M}\right)^2\\\ppsi_i\1=1~(i=1,\ldots, N)\\\ww_k\1=1/\Delta t~(k=1,\ldots, q)}}\OO\left(\PPsi,\WW\right)=\inf\limits_{\substack{\ww_1\in\xR_+^M\\\ww_1\1=1/\Delta t}}\cdots\inf\limits_{\substack{\ww_q\in\xR_+^M\\\ww_q\1=1/\Delta t}}\inf\limits_{\substack{\ppsi_1\in\xR_+^q\\\ppsi_1\1=1}}\cdots\inf\limits_{\substack{\ppsi_N\in\xR_+^q\\\ppsi_N\1=1}}\OO\left(\PPsi,\WW\right).
\end{equation*}
The idea is to minimize the criterion $\OO\left(\PPsi,\WW\right)$ working on one line-vector at a time: $\ppsi_1$, then $\ppsi_2$, and so on until $\ppsi_N$, and then $\ww_1$, $\ww_2$ and so forth until $\ww_q$. This process being repeated as many times as necessary -- and affordable -- to get some kind of convergence.

Let us make one more remark before giving out the algorithm. Let $i\in\{1,\ldots, N\}$. The two following programs are equivalent -- in the sense that their feasible and optimal solutions are the same (we just reduced the criterion to minimize).
\begin{equation*}
	\minimize{\OO\left(\PPsi,\WW\right)}{\ppsi_i}{\begin{cases}\ppsi_i\in\xR_+^q\\\ppsi_i\1=1\end{cases}}
	\enskip\Longleftrightarrow\enskip
	\minimize{\Frac{1}{2}\norm{\ff_i-\ppsi_i\WW}_{\xR^M}^2}{\ppsi_i}{\begin{cases}\ppsi_i\in\xR_+^q\\\ppsi_i\1=1 .\end{cases}}
\end{equation*}
The algorithm we implemented to numerically solve \eqref{vectorial minimization} is the following.
\begin{lgrthm}[AQM]\
	\begin{itemize}
		\item Initialize $\PPsi$ and $\WW$ with uniform values.
			\begin{align*}
				\psi_{ik}&=1/q&&(i=1,\ldots, N)(k=1,\ldots, q)
				\\\w_{kj}&=1/(M\Delta t)&&(k=1,\ldots, q)(j=1,\ldots, M).
			\end{align*}
		\item While some stopping criterion is not reached (\text{i.e.} a maximum number $iter_{\text{max}}\approx20$ of iterations),
			\begin{enumerate}[(a)]
				\item do, for $i = 1,\ldots, N$,
					\begin{equation*}
						\minimize{\Frac{1}{2}\norm{\ff_i-\ppsi_i\WW}_{\xR^M}^2}{\ppsi_i}{\begin{cases}\ppsi_i\in\xR_+^q\\\ppsi_i\1=1.\end{cases}}
					\end{equation*}
				\item do, for $k = 1,\ldots, q$,
					\begin{equation*}
						\minimize{\OO\left(\PPsi,\WW\right)}{\ww_k}{\substack{\begin{cases}\ww_k\in\xR_+^M\\\ww_k\1=1/\Delta t.\end{cases}}}
					\end{equation*}
			\end{enumerate}
		\end{itemize}
\end{lgrthm}
Each of the minimization problems addressed in this algorithm has a convex quadratic criterion depending on at most $M\approx2000$ variables constrained by the same amount of lower bounds plus one linear equality.
These optimization problems were solved with a dual method of the \emph{active set} type, detailed in Goldfarb and Idnani \cite{Goldfarb1983}, specifically designed to address strictly convex quadratic programs.
For a given quadratic program $P$, the \emph{active set} of a point $\xx$ is the set of (linear) inequality constraints satisfied with equality (\emph{active constraints}).
The algorithm starts from a point $\xx$ such that the set $A$ of its active constraints (possibly empty) is linearly independent, and the following steps are repeated until all constraints are satisfied.
First, a violated constraint $\mathbf{c}$ is picked.
Second, the feasibility of the subproblem $P(A\cup\{\mathbf{c}\})$ of $P$ constrained only by the constraint $\mathbf{c}$ and those of the current active set $A$ is tested.
If $P(A\cup\{\mathbf{c}\})$ is infeasible then $P$ is infeasible as well.
Otherwise $\xx$ is updated with a point whose active set is a linearly independent subset of $A\cup\{\mathbf{c}\}$ including $\mathbf{c}$.
The loop terminates with an optimal solution.



		\subsubsection{Metamodelling of functional decomposition coefficients}
			\label{section:coefficients}
			In Sections \ref{section:CPCA}, \ref{section:EIM} and \ref{section:AQM}, three methods have been studied to approximate a sample of probability density functions on basis such that:
\begin{equation*}
	\hat{f}_i=\sum_{k=1}^q\psi_{ik}w_k.
\end{equation*}
The functions are characterized by their coefficients. The problem dimension is therefore reduced to the basis size $q$. A metamodel must now be designed to link the input parameters and the coefficients: a function $\hat{\ppsi}:\X\to\xR^q$, where $\hat{\ppsi}=(\hat{\psi}_1,\ldots,\hat{\psi}_q)$. The metamodel is built thanks to the known points: $\psi_{ik}=\psi_k(\xx_i)$ ($i=1,\ldots, N$,$k=1,\ldots, q$). For a new $\xx\in\X$, we look for the estimation $\hat f_\xx$ of the output $f_\xx$ of $G$:
\begin{equation*}
  \hat{f}_\xx=\sum_{k=1}^q\hat{\psi}_k(\xx)w_k.
\end{equation*}

For MMP and AQM methods, the searched coefficients must respect the following constraints: for all $\xx\in\X$,
\begin{equation*}
	\left\{\begin{aligned}\sum\limits_{k=1}^q\psi_k(\xx)&=1\\\psi_k(\xx)&\geq0\qquad(k=1,\ldots,q).\end{aligned}\right.
\end{equation*}
Therefore, the outputs of the metamodel to be built must be positive and have their sum equal to one.

So far, the main approach investigated has been to build separate metamodels on each coefficient without constraining the output of the metamodel. When the metamodel is used to make predictions on new points in the input space, the predicted coefficients which are lower than zero are put equal to zero, and the coefficients are then renormalized such that their sum is equal to one. Any type of metamodel can be used in this approach. Polynomial regression, generalized additive models \cite{Hastie1990} and Gaussian process metamodels have been tested, but the obtained results on different test cases are not convincing. This approach on the presented test cases does not lead to good results. The metamodel error is bigger than the error due to the decomposition on a functional basis.

Future research should be dedicated on the construction of more efficient metamodels taking into account these constraints. Metamodels which ensure the positivity of the output exist in the literature (see for example \cite{Daveiga2012} for Gaussian process with inequality constraints). The main difficulty is to ensure that the sum of the coefficients is equal to one. A solution can be to consider the coefficients as realizations of a random variable on a simplex space. Aitchison \cite{Aitchison1982} proposes to apply a bijective transformation from this simplex space to $\mathbb R^q$. A metamodel without constraint could then be built between the input parameters and the transformed coefficients.


\section{Numerical tests on toy functions}
	\label{chapter:numerical_tests}
	In this section, we apply the methods presented on two toy examples respectively for $d=1$ in Section \ref{section:toy_example_1} and $d=5$ in Section \ref{section:toy_example_2}. We compare the relative error of the estimation of different quantities of interest: $\xLone$, $\xLn{2}$ and the Hellinger distances between two densities, mean, variance and the ($1\%,\ 25\%,\ 75\%$,\ $99\%$)-quantiles. We recall that the Hellinger distance between $f$ and $g$ is the $\xLn{2}$ distance between $\sqrt f$ and $\sqrt g$. The relative error is defined below for the three norms and the other scalar quantities of interest:

\begin{align*}
	100&\Frac{\int\limits_I\abs{f(t)-\hat{f}(t)}^2\xdif{t}}{\int_If(t)^2\xdif{t}}\\
	100&\Frac{\int\limits_I\abs{\sqrt{f(t)}-\sqrt{\hat{f}(t)}}^2\xdif{t}}{\int_I\sqrt{f(t)}^2\xdif{t}}\\
	100&\Frac{\int\limits_I\abs{f(t)-\hat{f}(t)}\xdif{t}}{\int_I\abs{f(t)}\xdif{t}}\\
	100&\Frac{\abs{u-\hat{u}}}{\abs{u}},
\end{align*}
where $u$ is the mean, variance and the studied quantiles.
This study is separated in two independent parts for the two toy examples.

First, we compare the error obtained from the two estimators $\hat{f}_{\xLn{2}}$ and $\hat{f}_\mathrm{He}$ based on kernel regression. We show the relative error in function of different sizes $N_1 < \cdots < N_k$ of the design of experiments such that a design with $N_i$ points is included in the designs with $N_j$ points $(i<j)$. For each $N_i$, relative error for different quantities of interest are averaged on the same 1000 test points chosen uniformly in $\X$. These computations are repeated 25 times with different designs of experiments. 

In the second part, we compare the relative error obtained from the construction of basis obtained in Sections \ref{section:CPCA}, \ref{section:EIM} and \ref{section:AQM} to reconstruct the $N$ density probability functions of the design of experiments. These two methods decompose the known probability density functions on a functional basis and approximate them. However, no metamodel has been designed between the coefficients of the probability density functions on the basis and the parameters of the computer code, so that no estimation of an unknown probability density function can be done. Therefore, both methods cannot be compared to the kernel regression.

Each probability density function is discretized on 512 points. 

	\subsection{A one-dimensional example (Toy example 1)}
		\label{section:toy_example_1}
		The first toy example is defined as follows:
\begin{eqnarray*}
    G(x, \xi_1,\xi_2,U) = (\sin(x  (\xi_1 + \xi_2)) + U)\mathds{1}_{ \{\sin(x  (\xi_1 + \xi_2)) + U \geq -1 \} }
\end{eqnarray*}
where
\begin{eqnarray*}
x & \in & \X = [0 , 1] \\
\xi_1 & \sim & \mathcal{N}(1,1) \\
\xi_2 & \sim & \mathcal{N}(2,1) \\
U & \sim & \mathcal{U}([0,1]).
\end{eqnarray*}

Let $x_0 \in [0, 1]$, we want to estimate the probability density function $f_0$ of the random variable $G(x_0, \xi_1,\xi_2,U)$. 

\begin{figure}
\center
		\includegraphics[scale = 0.35, bb=0 0 1076 587]{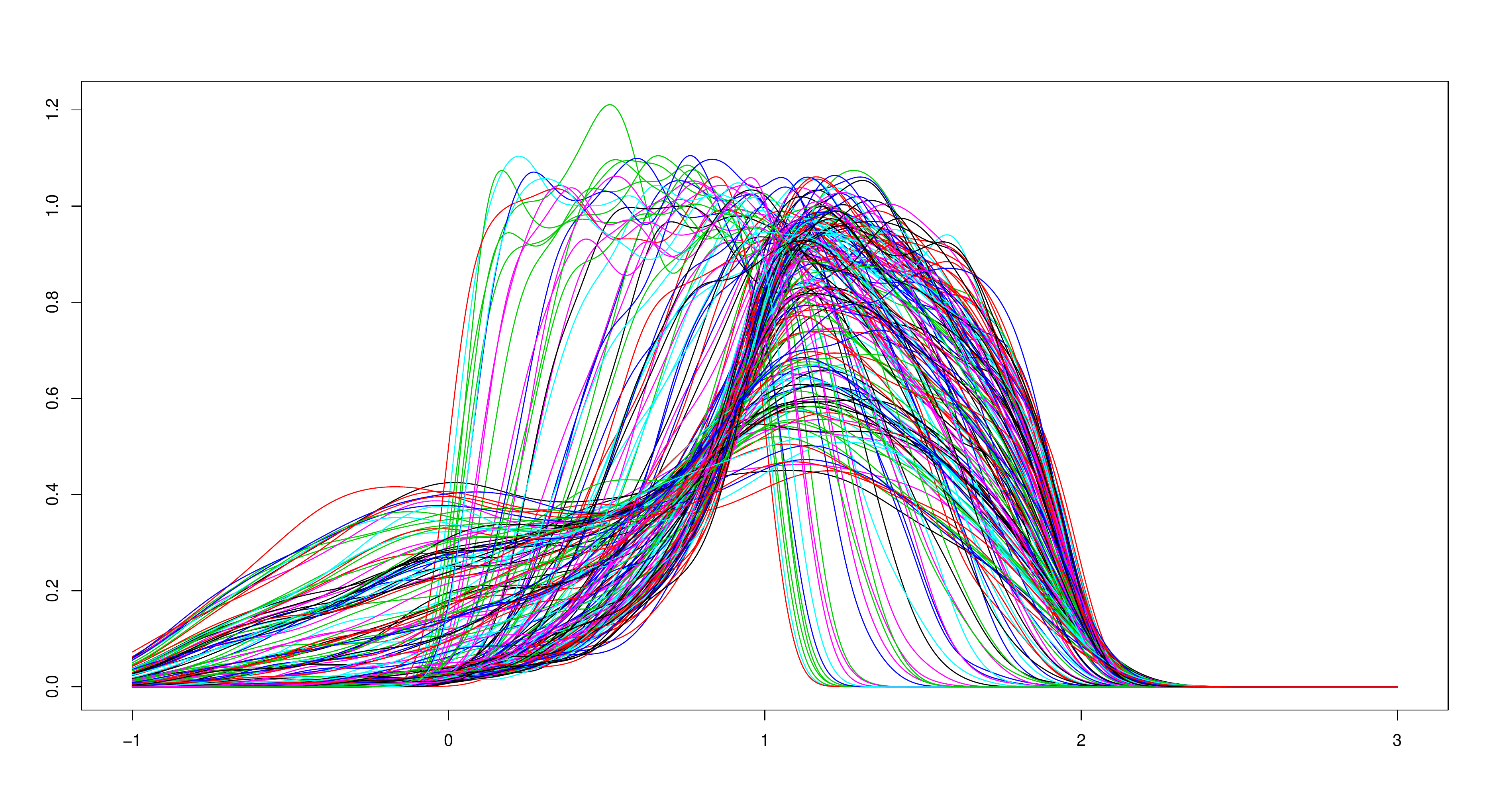}
	\caption{Toy example 1: representations of outputs for $N = 100$ different parameters.}
\label{figure:exemple:1:data}
\end{figure}

\subsubsection{Kernel regression}

\begin{figure}
\center
		\includegraphics[scale= 0.4, bb=0 0 720 720]{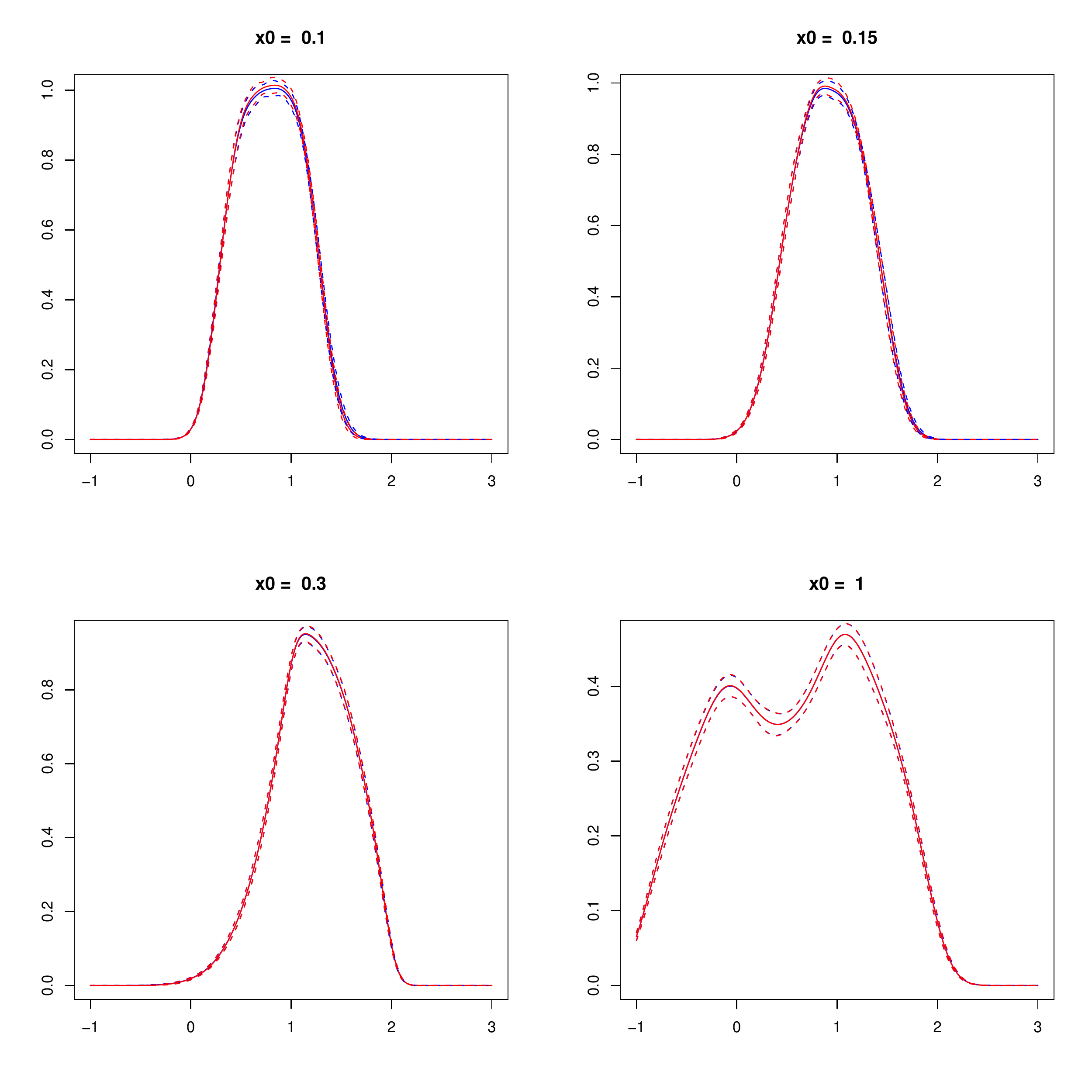}
	\caption{Toy example 1: curves in red (resp. blue) plain line represents the estimation of a density for a fixed $x_0$ from the first (resp. second) method and dashed red (resp. blue) line represent the standard deviation of the estimation obtained from 25 independent design of experiments. The blue and the red lines are superposed.}
\label{figure:exemple:1:1}
\end{figure}

\begin{figure}
\center
		\includegraphics[scale=0.4, bb=0 0 720 720]{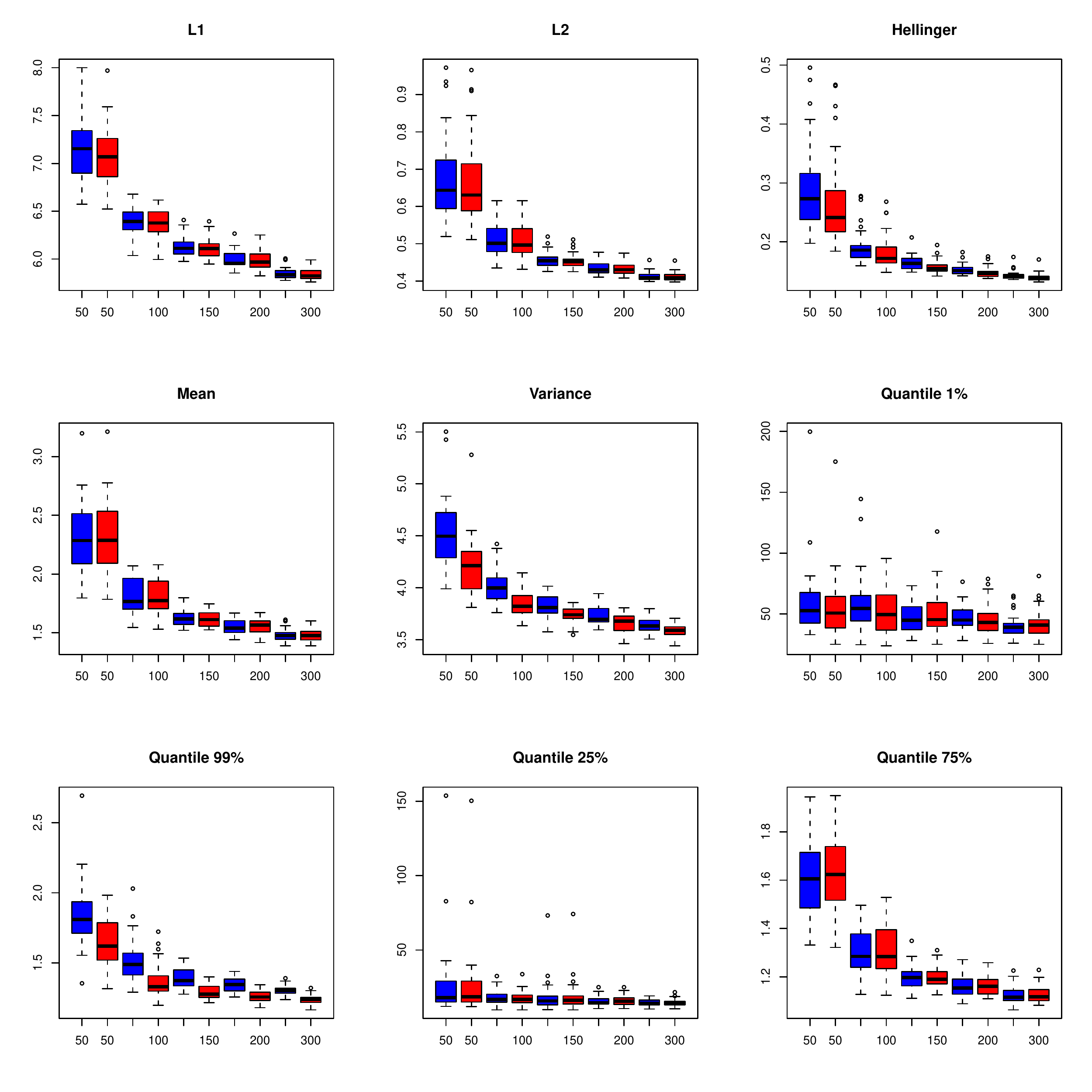}
	\caption{Toy example 1: boxplot of the errors for different sizes of $N$. Blue: estimator given by $\xLn{2}$ norm. Red: estimator given by the Hellinger distance. Results have been averaged with 25 independent experiments. Blue: estimator given by $\xLn{2}$ norm. Red: estimator given by the Hellinger distance.}
\label{figure:exemple:1:boxplot}
\end{figure}

In this section, one compares the estimators based on Hellinger distance and $\xLn{2}$ norm given respectively by equations (\ref{eq:L2:1}) and (\ref{eq:Hellinger:1}). Figure \ref{figure:exemple:1:data} represents the output for $N = 100$ points simulated uniformly and independently in $\X$. Figure \ref{figure:exemple:1:1} represents for different values of $x_0$ the mean (in plain line) and the standard deviation (in dashed line) of the estimation obtained from the two estimators of the kernel regression in red for $\xLn{2}$ norms and in blue for the Hellinger distance. For these four parameters, the two estimators give approximately the same estimations. Figure \ref{figure:exemple:1:boxplot} shows the boxplot of the relative error for $\xLn{2}$ norms and the Hellinger distance. In this first example the estimator based on the Hellinger distance seems slightly better than the one based on the $\xLn{2}$ norm. The errors are very close for the norms and the mean but are more different for other quantities of interest. The variance of the error decreases very quickly for norms. Moreover, the errors are low for most quantities of interest except for 1\% and 25\% quantiles. 

%

\subsubsection{Functional decomposition methods}

\begin{figure}
\center
		\includegraphics[scale=0.35, bb=0 0 1296 864]{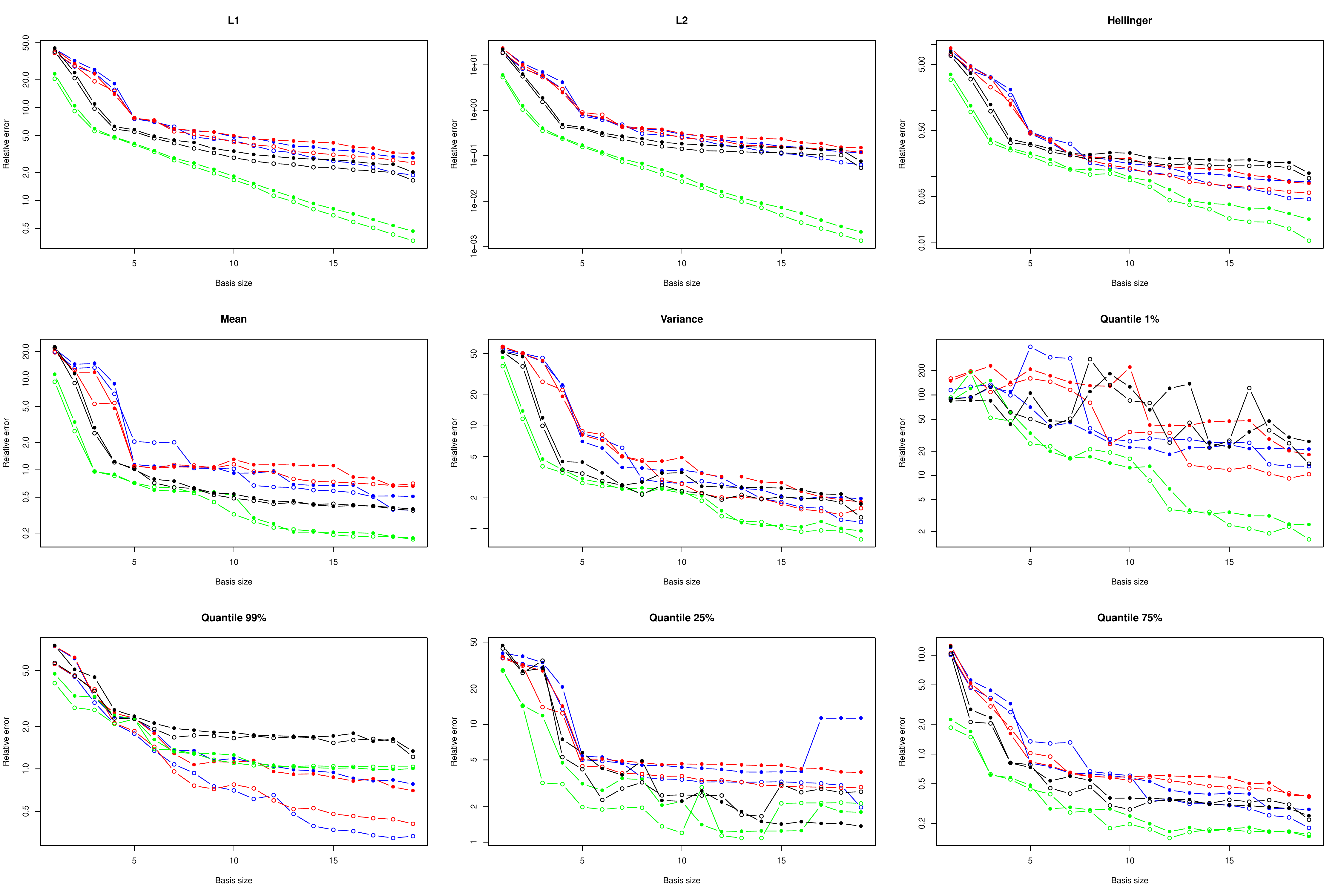}
	\caption{Toy example 1: comparison of the relative error for different quantities in function of the size of the basis $q$ with a design of experiments of size $N = 50$ (circles) and 100 (filled circles). Blue: MMP decomposition given by $\xLn{2}$ norm. Red: MMP decomposition given by the Hellinger distance. Black: AQM decomposition. Green: CPCA decomposition.}
\label{figure:exemple:1:decomposition}
\end{figure}

In this section, the four functional decompositions presented in \ref{section:DR}, MMP with $\xLn{2}$ and Hellinger distances, AQM and CPCA are compared. The decompositions are applied to learning samples whose sizes are $N=50$ and 100 and for decomposition basis sizes ranging from 1 to 20. The relative errors between the learning sample functions and their approximations are computed for the 9 quantities of interest. Figure \ref{figure:exemple:1:decomposition} represents, in logarithmic scale, these relative errors in function of the basis size for the four decompositions. First, the relative errors are low for all quantities of interest and especially for the norms, the mean and the higher quantiles. For the four decompositions, the relative error for norm decreases very evenly. The decrease is quite regular too for the mean, variance, 75\% and 99\% quantiles. The quality of approximation of the lower quantiles behaves more irregularly, as it sometimes increases with the basis size. CPCA method outperforms the three others for the three distances. This method seems better for all other quantities of interest except the 25\% quantile, but the gap between the decompositions is smaller.
The relative error for the learning sample of 100 densities is for most of the quantities higher than the relative error for the sample of 50 densities. For the same number of components $q$, it is indeed more difficult to approximate a higher number $N$ of functions, so that the relative error increases with the size of the design of experiments for a constant number of components.

\begin{figure}
\center
		\includegraphics[scale = 0.4, bb=0 0 576 360]{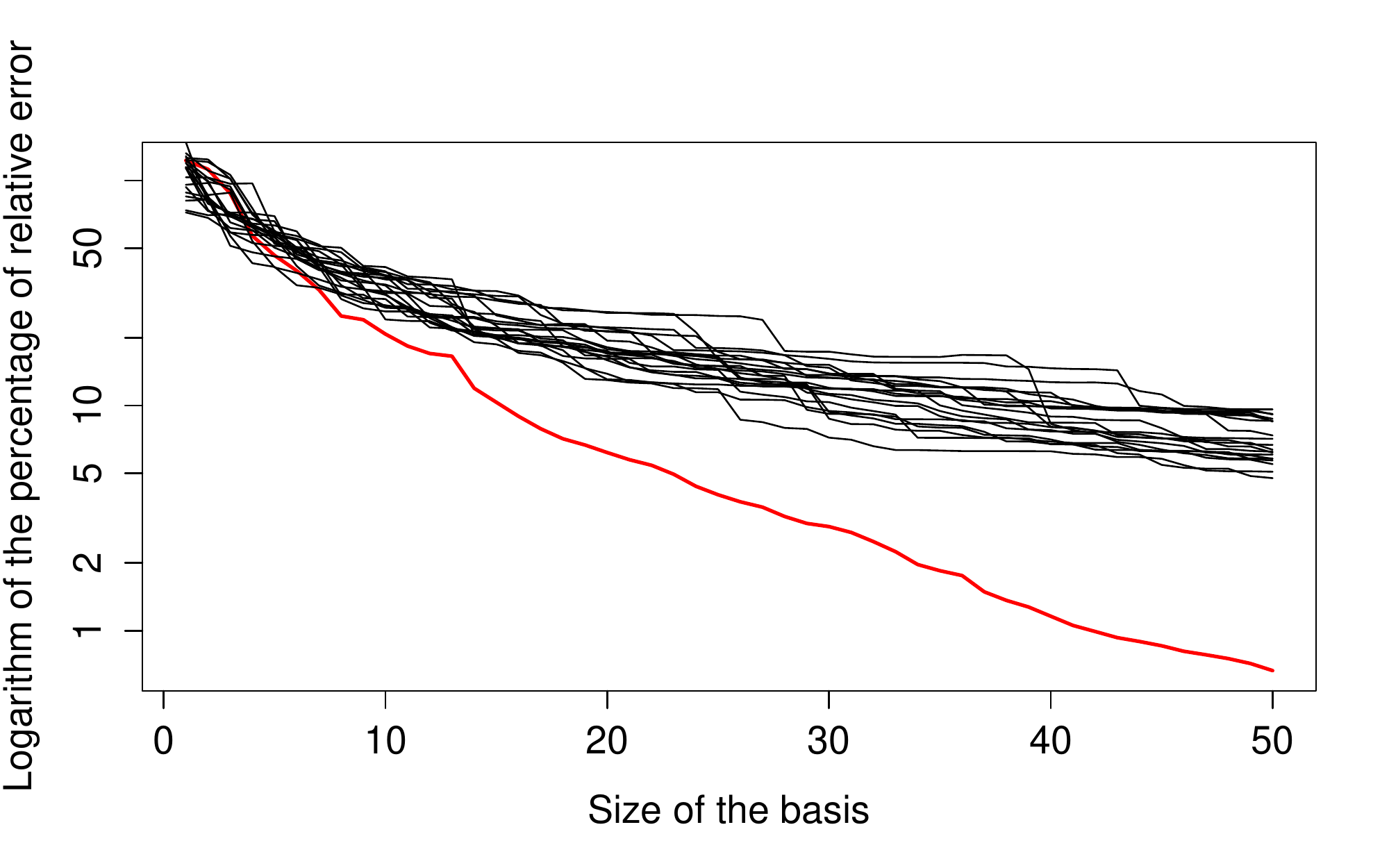}
	\caption{Toy example 1: comparison of the $\xLn{2}$ relative error obtained with MMP (in red) and a random choice in the design of experiments (in black).}
	\label{figure:exemple:1:EIMvsAlea}
\end{figure}

Figure \ref{figure:exemple:1:EIMvsAlea} represents the $\xLn{2}$ relative error of MMP (in red) with $N = 200$ and the relative error of twenty independent basis (in black) containing probability functions chosen randomly in the learning sample $f_1,\ldots,f_N$. The MMP outperforms the random strategies. This result validates the way to choose a new curve among the sample of probability density functions, as this choice is better than randomly picking functions in the sample.

	\subsection{A five-dimensional example (Toy example 2)}
		\label{section:toy_example_2}
\begin{figure}
\center
		\includegraphics[scale=0.4, bb=0 0 1076 587]{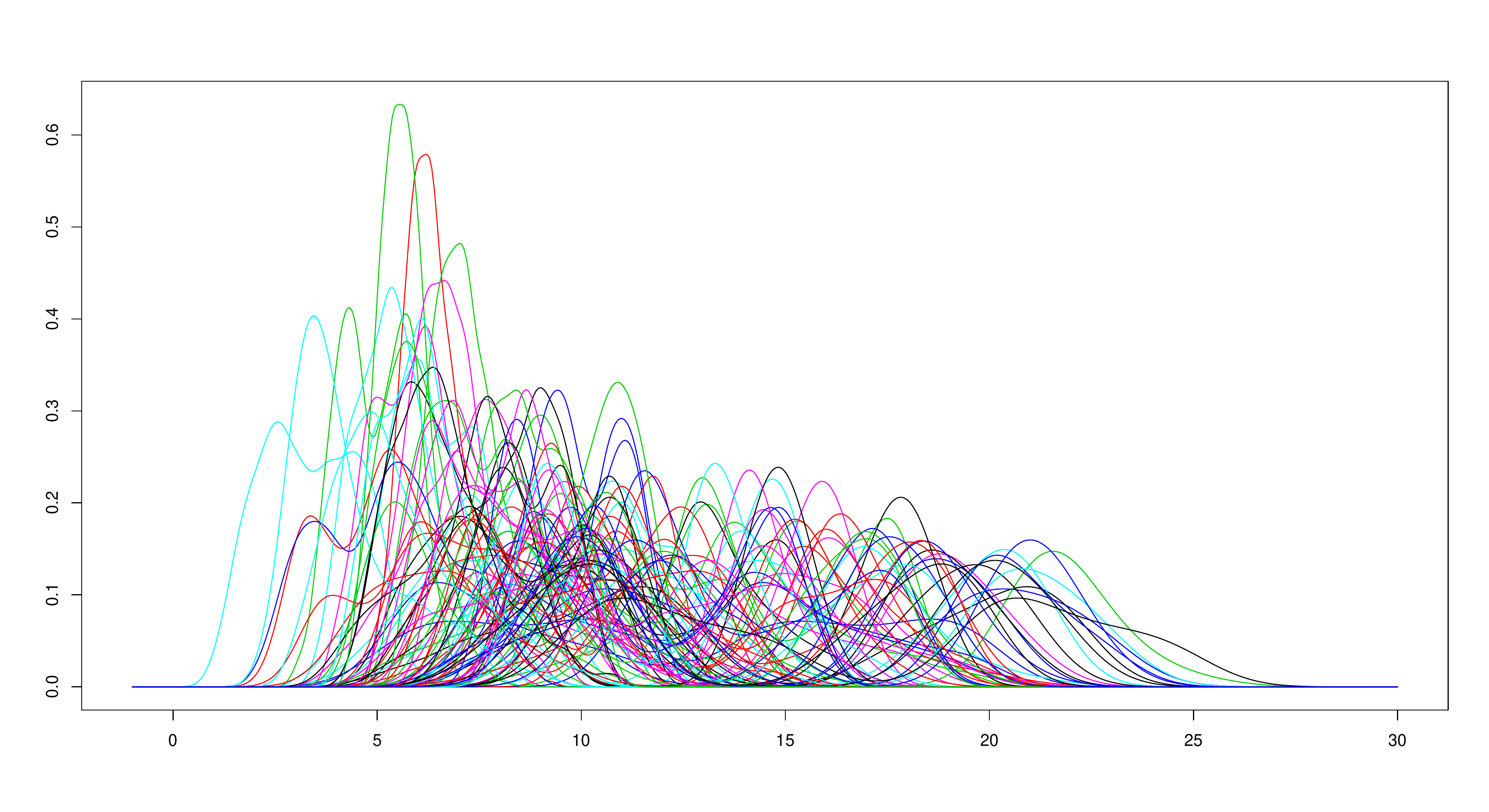}
	\caption{Toy example 2: representations of outputs for $N = 100$ different parameters.}
\label{figure:exemple:2:data}
\end{figure}

Let $\xx = (x_1,x_2,x_3,x_4,x_5)\in\X = [0 , 1]^5$. A second toy example is defined as follows:
\begin{eqnarray*}
    G(\xx,N,U_1,U_2,B) = (x_1+ 2x_2 + U_1)\sin(3x_3- 4x_4+N) + U_2+10x_5B + \sum\limits_{i=1}^5 ix_i
\end{eqnarray*}
where
\begin{eqnarray*}
N & \sim & \mathcal{N}(0,1) \\
U_1 & \sim & \mathcal{U}([0,1]) \\
U_2 & \sim & \mathcal{U}([1,2]) \\
B & \sim & \mathcal{B}ern(1/2).
\end{eqnarray*}
Figure \ref{figure:exemple:2:data} represents the output for $N = 100$ points simulated uniformly and independently in \X.

\subsubsection{Kernel regression}

\begin{figure}
\center
		\includegraphics[scale=0.4, bb=0 0 1076 587]{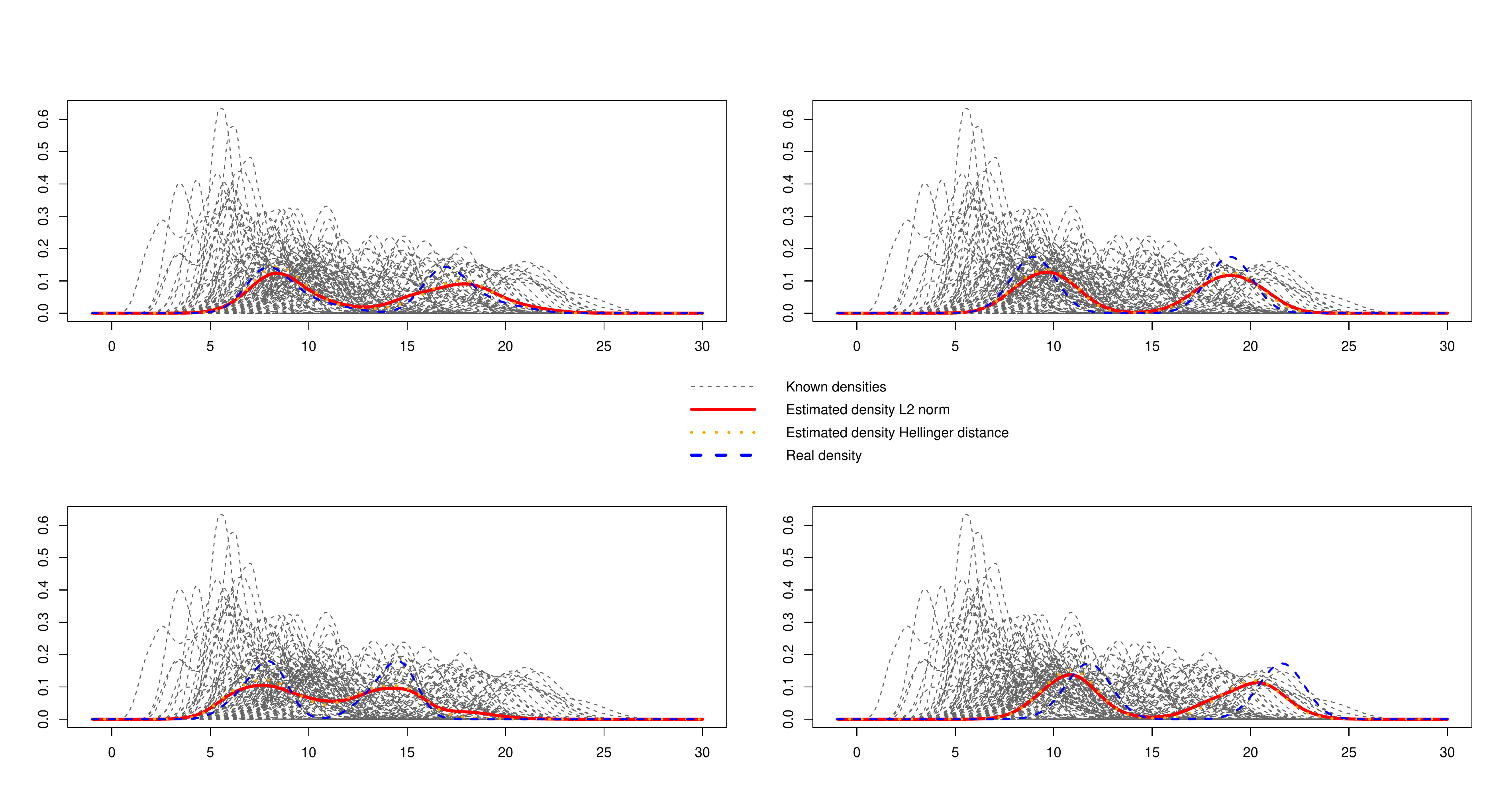}
	\caption{Toy example 2: curves in dashed gray lines are probability density functions, red plain line (resp. orange dots line) represents the estimation of a density for a fixed $\xx_0$ from the first (resp. second) method and dashed blue line is the true probability density function. These estimations were made with $h$ isotropic. The orange dotted line and the red line are superposed.}
\label{figure:exemple:2:1}
\end{figure}

\begin{figure}
\center
		\includegraphics[scale = 0.4, bb=0 0 720 720]{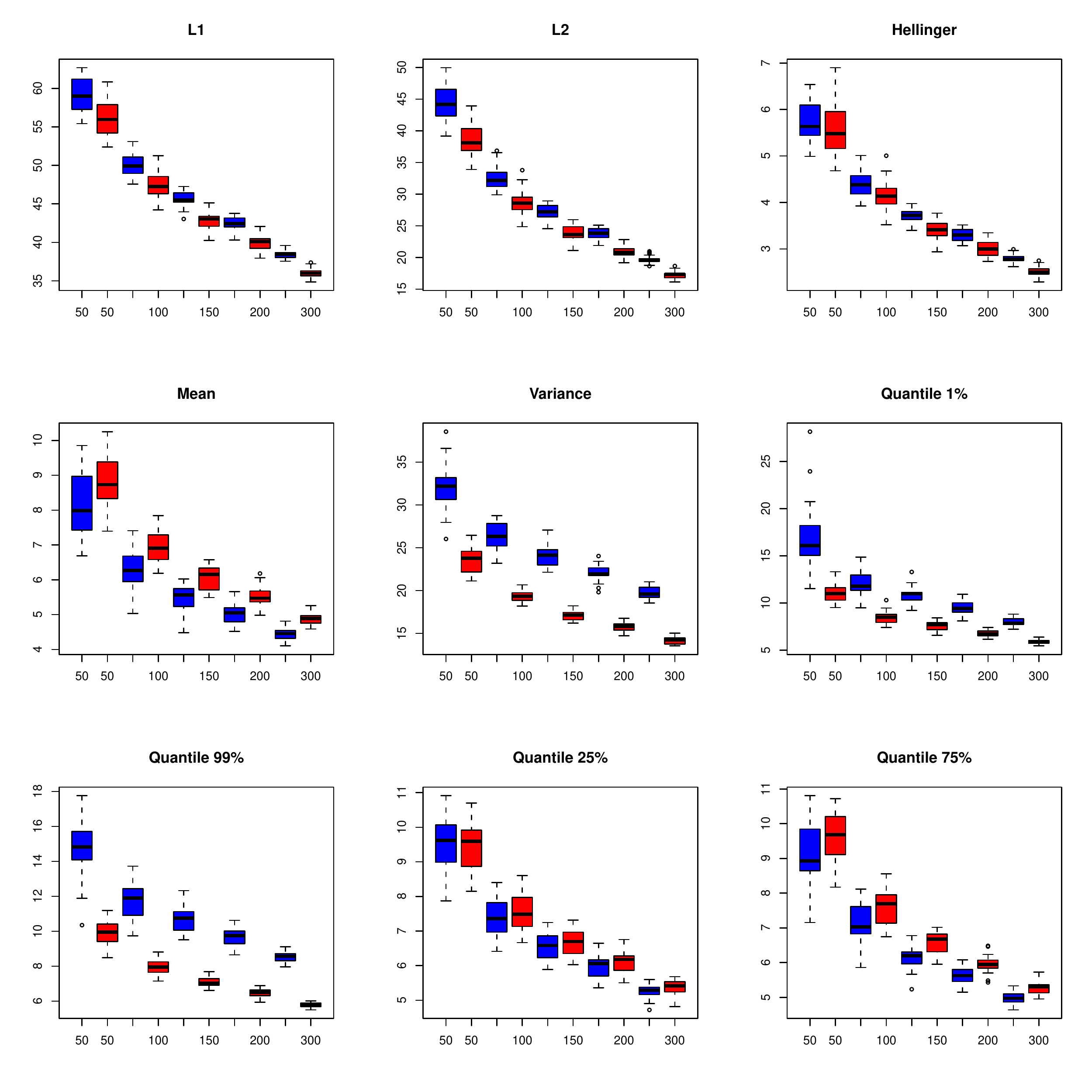}
	\caption{Toy example 2: boxplot of the errors for different sizes of $N$. Left: estimator given by $\xLn{2}$ norm. Right: estimator given by the Hellinger distance. Results have been averaged with 25 independent experiments. These estimations were made with $h$ isotropic. Blue: estimator given by $\xLn{2}$ norm. Red: estimator given by the Hellinger distance.}
\label{figure:exemple:2:boxplot:uni}
\end{figure}

\begin{figure}
\center
		\includegraphics[scale = 0.4, bb=0 0 720 720]{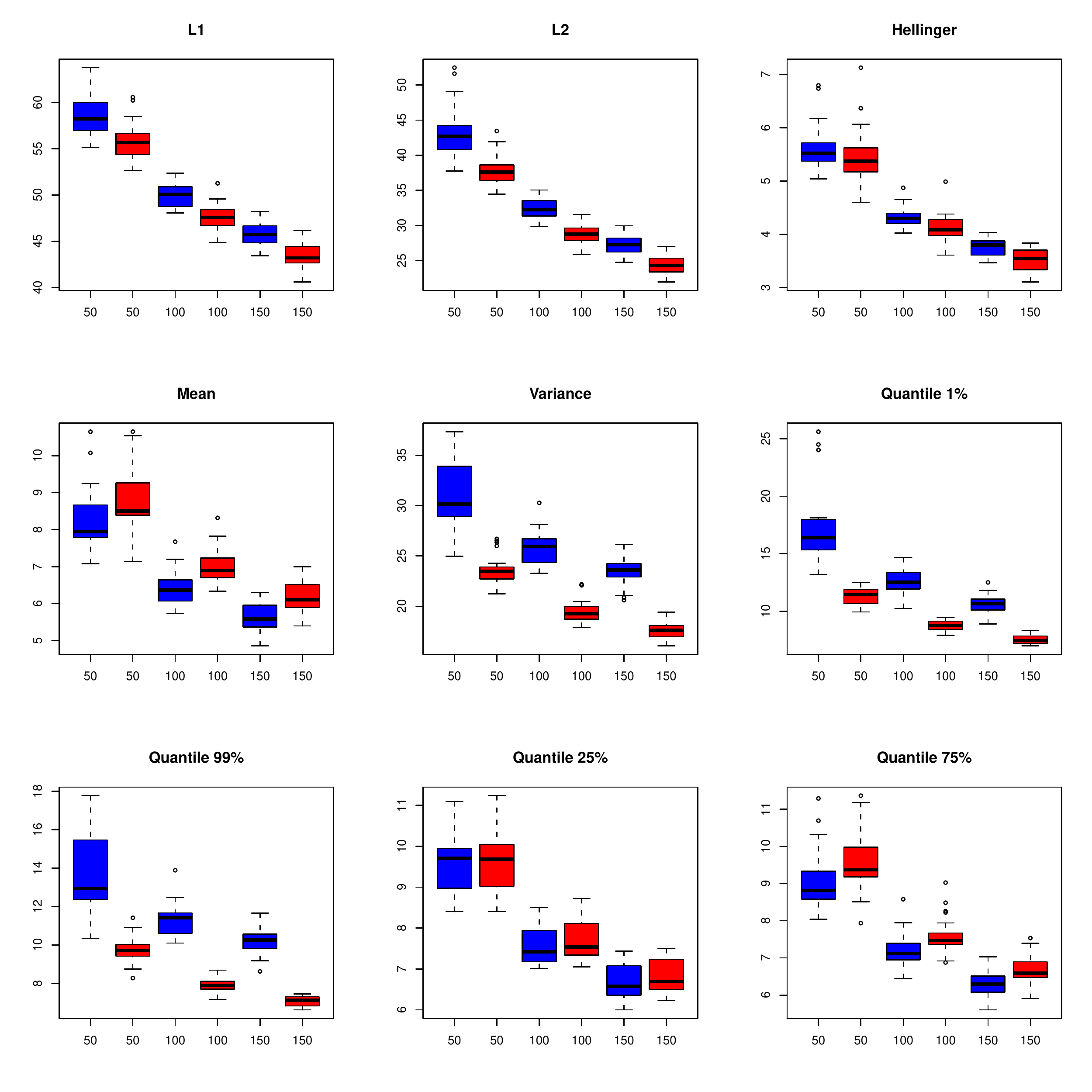}
	\caption{Toy example 2: boxplot of the errors for different sizes of $N$. Left: estimator given by $\xLn{2}$ norm. Right: estimator given by the Hellinger distance. Results have been averaged with 25 independent experiments. These estimations were made with $h$ anisotropic. Blue: estimator given by $\xLn{2}$ norm. Red: estimator given by the Hellinger distance.}
\label{figure:exemple:2:boxplot:multi}
\end{figure}

In this section, the  kernel regression is tested on this second toy example with isotropic (Figure \ref{figure:exemple:2:boxplot:uni}) and anisotropic bandwidth (Figure \ref{figure:exemple:2:boxplot:multi}). As in the first example, Figure \ref{figure:exemple:2:1} represents for different values of $x_0$ the true probability density function (blue dashed line) and the estimation obtained from the two estimators of the kernel regression (red plain line and orange dots line). In this example, the two estimators give different results.
In red is represented the relative error for the estimator based on Hellinger distance, in blue for the estimator based on $\xLn{2}$ norm. The Hellinger estimator gives a better approximation of the quantities presented in Figure \ref{figure:exemple:2:boxplot:uni}, except for the mean. Contrary to the case of the first toy example, the errors are quite high.
The variance of the error on the norms, presented in Figure \ref{figure:exemple:2:boxplot:uni}, is low but does not seem to decrease steadily. 
The use of an anisotropic bandwidth does not seem to improve the quality of the estimation, while it is much longer to compute. The errors are approximately the same with isotropic or anisotropic bandwidth. The greater precision brought by the anisotropic bandwidth is compensated by the difficulty to estimate it.

%

%

\subsubsection{Functional decomposition methods}

\begin{figure}
\center
		\includegraphics[scale=0.35, bb=0 0 1296 864]{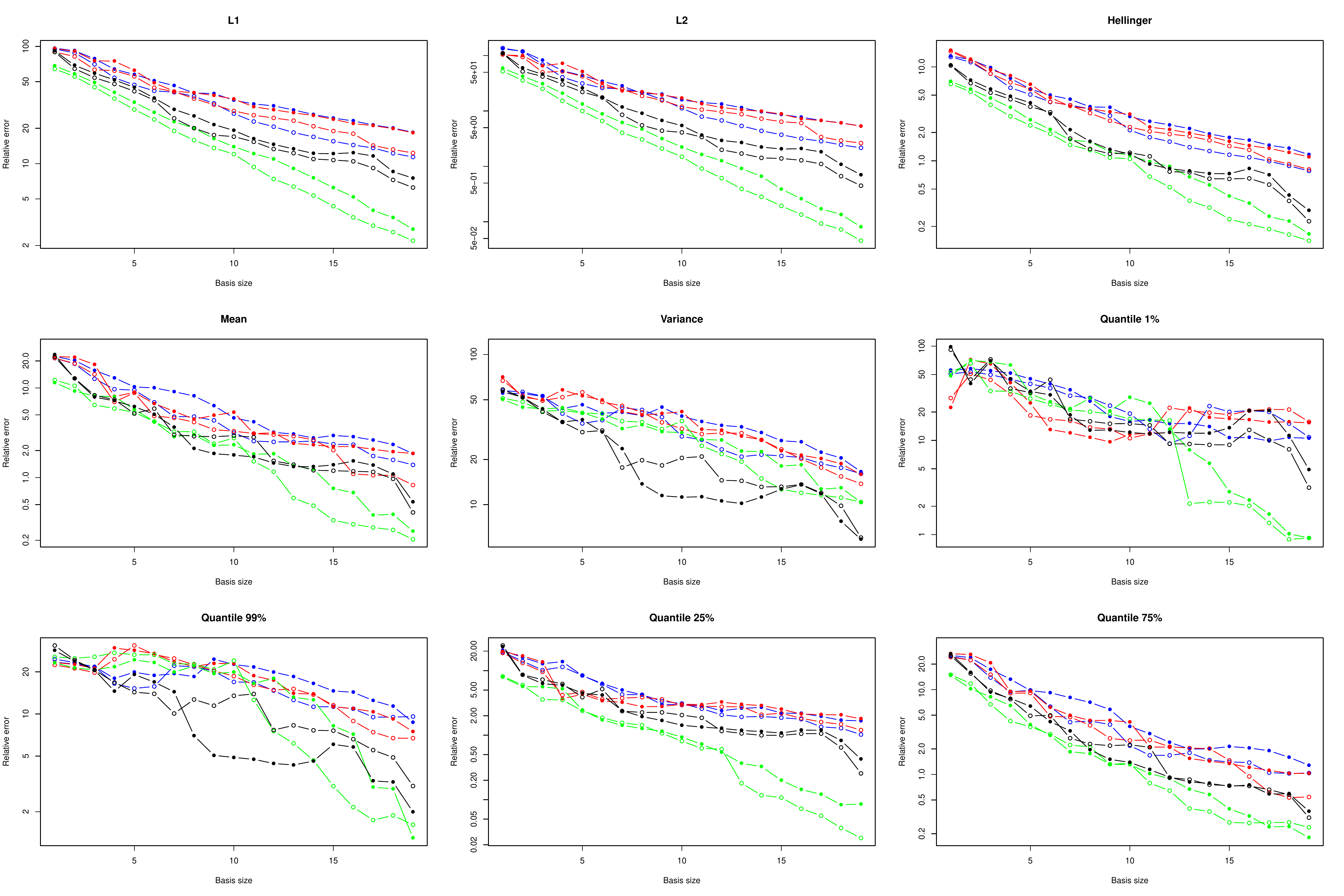}
	\caption{Toy example 2: comparison of the relative error for different quantities in function of the size of the basis $q$ with a design of experiments of size $N = 50$ (circles) and 100 (filled circles). Blue: MMP decomposition given by $\xLn{2}$ norm. Red: MMP decomposition given by the Hellinger distance. Black: AQM decomposition. Green: CPCA decomposition.}
\label{figure:exemple:2:decomposition}
\end{figure}

In this section, the four functional decompositions MMP with $\xLn{2}$ and Hellinger distances, AQM and CPCA are compared. The decompositions are applied to learning samples whose sizes are $N=50$ and 100 and for decomposition basis sizes ranging from 1 to 20. The relative errors between the learning sample functions and their approximations are computed for the 9 quantities of interest. Figure \ref{figure:exemple:2:decomposition} represents, in logarithmic scale, these relative errors in function of the basis size for the three decompositions. 
First, the relative errors are low for all quantities of interest, as in the previous example. 
The relative errors decrease is less even than in the one-dimensional case. The behavior of the errors on 1\% and 99\% quantiles is particularly unsteady. CPCA method gives better results than the other methods for $\xLn{1}$, $\xLn{2}$ norms, 25\% and 75\% quantiles. However, for all studied quantities, the AQM and CPCA method seem to give quite equivalent results. MMP method errors are slightly higher than AQM errors for most of the quantities.



\section{Industrial applications}
	\label{chapter:industrial_applications}
	In this section we apply the methods we developed in Section \ref{Chapter:theory} to three industrial numerical codes of CEA, EDF, and IFPEN. We compute the same relatives errors as for the toy examples of Section \ref{chapter:numerical_tests}.
	\subsection{CEA application: CASTEM test case}
		\label{section:cea}

\label{section:castem}
\subsubsection{Code description}

In the framework of nuclear plant risk assessment studies, the evaluation of component reliability during accidental conditions is a major issue required for the safety case. Thermal-hydraulic (T-H) and thermal-mechanic (T-M) system codes model the behaviour of the considered component subjected to highly hypothetic accidental conditions. In the study that we consider here, the T-M code CASTEM takes as input 13 uncertain parameters, related to the initial plant conditions or to the safety system characteristics. Three of them are functional T-H parameters which depend on time: fluid temperature, flow rate and pressure. The other ten parameters are T-M scalar variables. For each set of parameters, CASTEM calculates the absolute mechanical strength of the component and the thermo-mechanical actual applied load. From these two elements, a safety margin (SM) is deduced.

The objective is to assess how these uncertain parameters can affect the code forecasts and more specifically the predicted safety margin. However, CASTEM code is too time expensive to be directly used to conduct uncertainty propagation studies or global sensitivity analysis based on sampling methods. To avoid the problem of huge calculation time, it can be useful to replace CASTEM code by a metamodel.
One way to fit a metamodel on CASTEM could be to discretize the functional inputs and to consider the values of the discretization as scalar inputs of CASTEM code. Nevertheless, this solution is often intractable due to the high number of points in the discretization. To cope with this problem, in \cite{Iooss2009, Marrel2012} a method was proposed to treat implicitly these \og uncontrollable\fg{} parameters functional parameters, while the other ten scalar parameters are considered as \og controllable\fg. CASTEM output is then a random variable conditionally to \og controllable\fg{} parameters.

\begin{figure}
\centering
		\includegraphics[scale = 0.45, bb=0 0 1080 432]{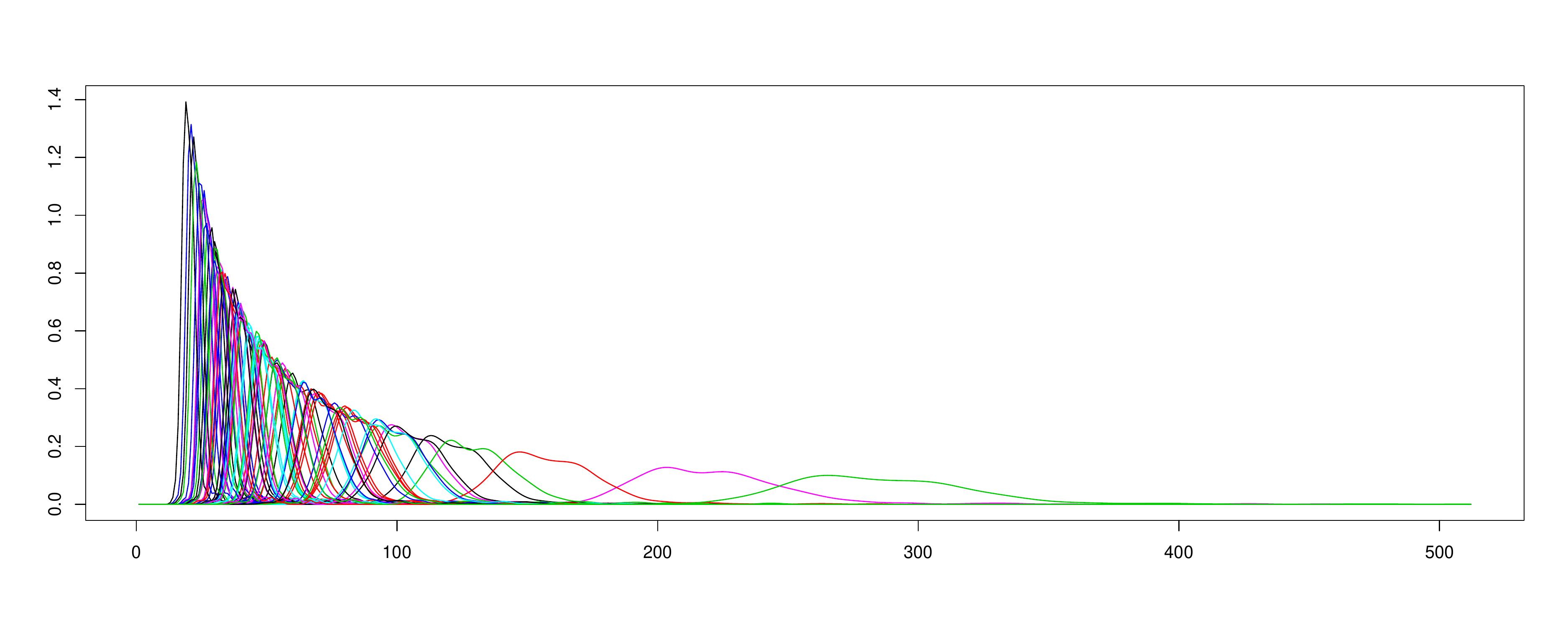}
	\caption{CASTEM test case: representations of outputs for $N = 75$ different parameters.}
\label{fig:castem:data}
\end{figure}

A latin hypercube sampling method \cite{McKay1979} is used to build a learning sample of 500 points in dimension 10. For each set of controllable parameters, CASTEM has been run 400 times with different uncontrollable functional parameters, randomly chosen in an available database. The probability density function $f_i$ ($i=1,..,500$) of the safety margin is computed by kernel estimation with the 400 outputs of CASTEM for each set of parameters. A few examples of the obtained probability density functions are represented on Figure \ref{fig:castem:data}. In the following, the two kernel regression metamodels, then MMP, AQM and CPCA decomposition methods are applied on CASTEM test case.

\subsubsection{Kernel regression}
\label{castem:kr}
\begin{figure}
\centering
		\includegraphics[scale = 0.45, bb=0 0 720 576]{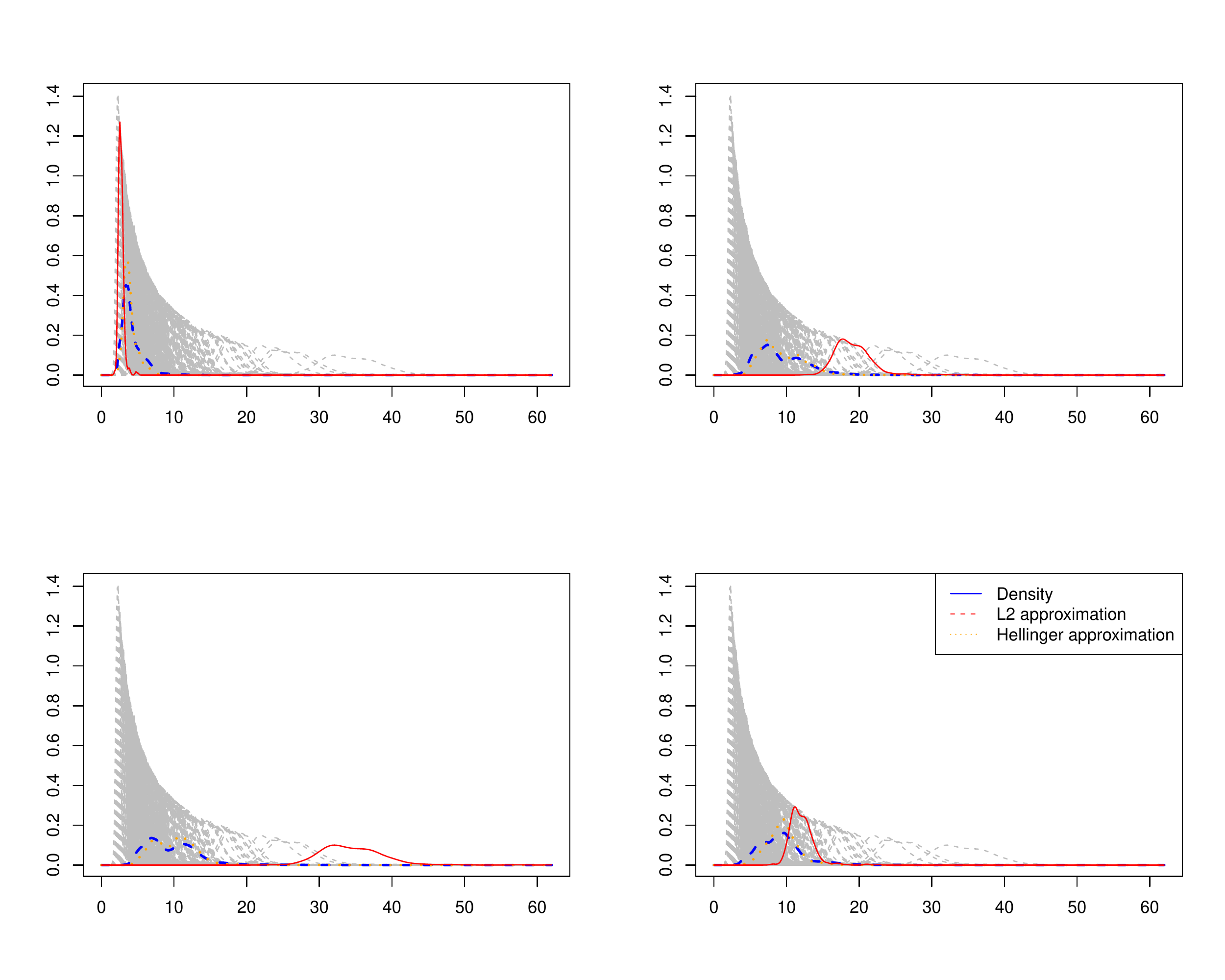}
	\caption{CASTEM test case: curves in dashed gray lines are probability density functions, red plain line represents the estimation of a density for a fixed $x_0$ and dashed blue line is the true probability density function.}
\label{fig:castem:4}
\end{figure}

In this section, the kernel regression method is applied with isotropic and anisotropic bandwidths. Figure \ref{fig:castem:4} represents for four different values of $x_0$ the true probability density function (blue dashed line) and the estimation obtained from the two estimators of the kernel regression ($\xLn{2}$ estimator in red plain line and Hellinger estimator in orange dotted line), with isotropic bandwidth. One can see that the estimation by the $\xLn{2}$ estimator of the four probability density functions is very far from the real function. In particular, at the bottom left of the figure, the mean of the predicted probability function in red is very far from the mean of the real one in dashed blue. On these four probability density functions, the Hellinger estimator gives much better results than the $\xLn{2}$ one.

To verify this first graphical analysis, the Leave-One-Out method has been used to assess the efficiency of the kernel regression-based metamodel. The bandwidth of the kernel regression is estimated thanks to all probability density functions except the $i^{\mathrm{th}}$ one. The function $f_i$ is estimated with the corresponding metamodel and the relative errors on the quantities of interest are computed between $f_i$ and its estimation. This process is repeated for each probability density function in the dataset and all the computed relative errors are averaged. 
The mean relative errors of the estimators based on Hellinger distance and $\xLn{2}$ norm are given in Table \ref{tab:castem} and Table \ref{tab:castem:anisotropic} for the different quantities of interest respectively for an isotropic and anisotropic bandwidths.
The two estimators perform poorly on this dataset. The errors for $\xLone$ and $\xLn{2}$ norms are particularly high. The Leave-One-Out validation confirms that the kernel regression method is not adapted to CASTEM test case. This can be explained by the important influence of the controllable T-M parameters compared to the one of the uncontrollable T-H parameters.

\begin{table}
\centering
\begin{tabular}{lcc}
\hline
 & Hellinger distance estimator & $\xLn{2}$ norm estimator\\
\hline
$\xLone$ norm        & 114.2\% & 114.9\%\\
$\xLn{2}$ norm       & 167.7\% & 171.9\%\\
Hellinger distance   & 17.0\% & 17.0\%\\
Mean                 & 22.2\% & 22.4\%\\
Variance             & 86.1\% & 86.9\% \\
1\% quantile         & 72.0\% & 74.3\%\\
99\% quantile        & 39.7\% & 40.9\%\\
25\% quantile        & 30.9\% & 31.5\%\\
75\% quantile        & 23.7\% & 23.9\%\\
\hline
\end{tabular}
\caption{CASTEM test case: the mean relative errors on the quantities of interest computed by the Leave-One-Out method for the kernel regression estimators based on the Hellinger distance and $\xLn{2}$ norm with an isotropic bandwidth.}
\label{tab:castem}
\end{table}

\begin{table}
\centering
\begin{tabular}{lcc}
\hline
 & Hellinger distance estimator & $\xLn{2}$ norm estimator\\
\hline
$\xLone$ norm        & 94.3\% & 94.5\%\\
$\xLn{2}$ norm       & 108.0\% & 108.6\%\\
Hellinger distance   & 14.3\% & 14.3\%\\
Mean                 & 15.4\% & 15.6\%\\
Variance             & 71.5\% & 71.9\% \\
1\% quantile         & 42.7\% & 43.3\%\\
99\% quantile        & 23.0\% & 23.8\%\\
25\% quantile        & 20.4\% & 20.4\%\\
75\% quantile        & 15.9\% & 16.3\%\\
\hline
\end{tabular}
\caption{CASTEM test case: the mean relative errors on the quantities of interest computed by the Leave-One-Out method for the kernel regression estimators based on the Hellinger distance and $\xLn{2}$ norm with an anisotropic bandwidth.}
\label{tab:castem:anisotropic}
\end{table}

\subsubsection{Functional decomposition methods}

\begin{figure}
\centering
		\includegraphics[scale = 0.35, bb=0 0 1296 864]{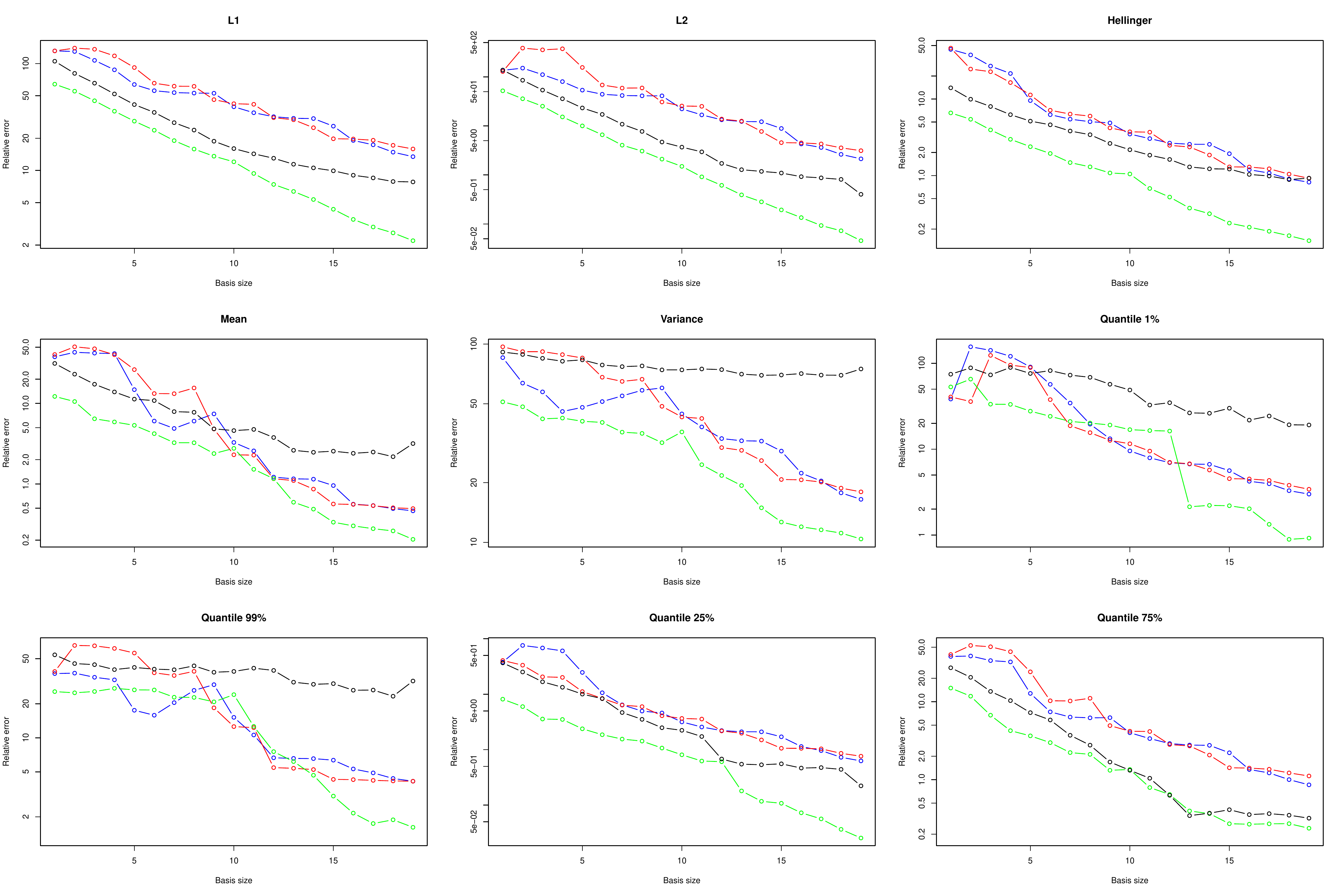}
	\caption{CASTEM test case: comparison of the relative error for the 9 quantities in function of the size of the basis $q$. Blue: MMP decomposition given by $\xLn{2}$ norm. Red: MMP decomposition given by the Hellinger distance. Black: AQM decomposition. Green: CPCA decomposition.}
\label{fig:castem:decompositions}
\end{figure}

In this section, the four functional decompositions MMP with $\xLn{2}$ and Hellinger distances, AQM and CPCA are compared. Figure \ref{fig:castem:decompositions} represents, in logarithmic scale, the relative errors on the different norms, modes and quantiles versus the basis size (from 1 to 20).
First, the relative errors are low for all quantities of interest, except for the variance. For the variance, the errors are over 10\% with a decomposition basis with 20 functions. The errors for small bases are high but decrease very quickly. The decrease is especially quick and even for errors on $\xLn{1}$, $\xLn{2}$ and Hellinger distances.
For modes and quantiles, the errors of the four methods do not decrease steadily for all quantities of interest, and especially for small basis sizes. 
CPCA clearly outperforms other methods for the $\xLn{1}$, $\xLn{2}$ and Hellinger distances. For other quantities of interest, it gives good results. For the variance, 1\% and 99\% quantiles, AQM method performs poorly compared to the three others. The errors of AQM and CPCA methods seem more stable than for MMP method for most of the quantities of interest, especially for modes and quantiles.
Overall, the results are quite good for the four methods. 

	\subsection{EDF application: VME test case}
		\label{section:edf}
		\subsubsection{Code description}
\begin{figure}
\center
		\includegraphics[scale = 0.45, bb=0 0 659 441]{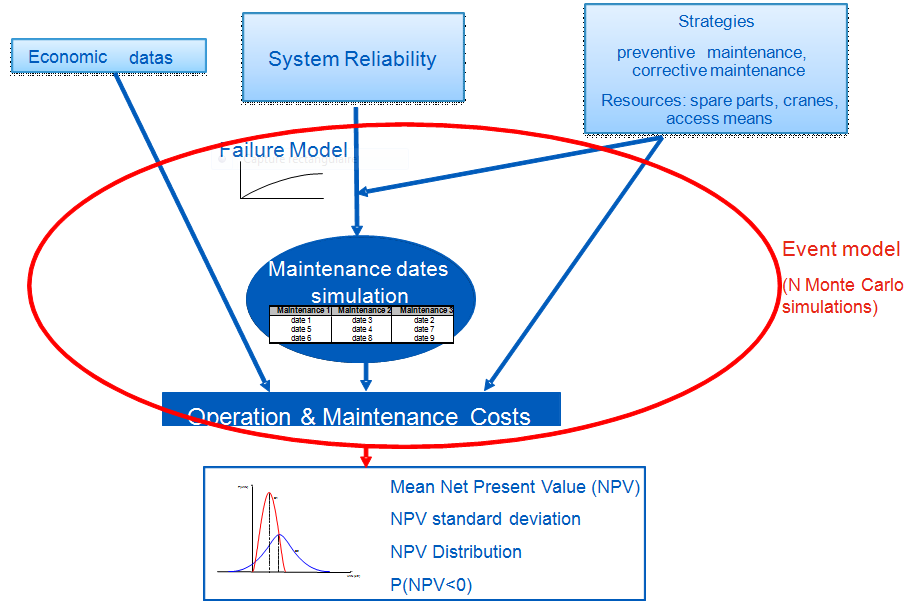}
	\caption{VME: EDF valorization of maintenance strategies tool description.}
\label{figure:VME:1}
\end{figure}

To optimize the whole life cost of its nuclear fleet, EDF has developed an asset management methodology \cite{Fessart2010}. A part of this methodology deals with exceptional maintenance tasks strategies. To help the decision maker to choose the best strategy (how many times do we need to carry out exceptional tasks, when,...?),  EDF has developed a dedicated tool called VME (described in Figure \ref{figure:VME:1}) based on Monte-Carlo simulation to compute many technical economic indicators among which the density function of the Net Present Value (called VAN for ``Valeur Actuelle Nette'') is the most relevant.

This tool leads to an important simulation time and requires an important amount of input data that are surrounded with uncertainties:
  \begin{itemize}
		\item Reliability data: generally there is not enough (or sometimes not any) feedback data to precisely evaluate reliability model parameters;
		\item Economic data: economic indicators and duration of maintenance tasks remain on several hypothesis that can be modified;
		\item Other data: uncertainty on operating times of power plants, maintenance tasks dates, etc. 
  \end{itemize}

In our work, we study a particular scenario simulation with VME concerning the life cycle management of four large transformers, adapted from the EPRI study made for CENG \cite{Fessart2011}. The uncertainty on the inputs are represented by four Weibull distributions which are summarized in Table \ref{tab:VME:1}. Define for $i=1,\ldots, 4$,
\begin{align*}
D_i = \{ \xx = (x_1,x_2,x_3) \in \mathbb{R}^3: \  \mathbb{E}[\mathcal{W}eibull(x_1,x_2, x_3 )] \in [ \mathbb{E}[W_i] \pm 0.2 \mathbb{E}[W_i]] \},
\end{align*}
and denote $sh_i, sc_i, loc_i$ respectively the shape, the scale and the location of $W_i$. For $i=1,\ldots,4$ let $x_{i,1},x_{i,2},x_{i,3},x_{i,4}$ be parameters simulated uniformly on $D_i$. Note that $W_2$ and $W_4$ follow the same law but $sh_2, sc_2, loc_2$ and $sh_4, sc_4, loc_4$ are considered as different parameters. That means there are three parameters for each Weibull distribution. Moreover, there are four other parameters which represent some costs and two other parameters which represent durations of services. All these $18$ parameters are summarized in Table \ref{tab:VME:2}.

\begin{table}
\center
\begin{tabular}{ccccc}
\hline
Uncertainty inputs  & distribution        & shape   & scale & location \\
$W_1$               & $\mathcal{W}eibull$ & 1/0.021 & 2.6   & -1  \\
$W_2$               & $\mathcal{W}eibull$ & 100     & 3.8   & -20 \\
$W_3$               & $\mathcal{W}eibull$ & 100     & 4     & -24 \\
$W_4$               & $\mathcal{W}eibull$ & 100     & 3.8   & -20 \\
\hline
\end{tabular}
\caption{VME: uncertainty inputs used in the model. Each distribution depends of three parameters: the shape, the scale and the location.}
\label{tab:VME:1}
\end{table}

\begin{table}
\center
\begin{tabular}{cccc}
\hline
Uncertainty inputs        & distribution        & support   \\
($sh_1, sc_1, loc_1$)     & $\mathcal{U}niform$ & $D_1$ &\\
($sh_2, sc_2, loc_2$)     & $\mathcal{U}niform$ & $D_2$ &\\
($sh_3, sc_3, loc_3$)     & $\mathcal{U}niform$ & $D_3$ &\\
($sh_4, sc_4, loc_4$)     & $\mathcal{U}niform$ & $D_4$ &\\
Cost1     & $\mathcal{U}niform$ & [796  ,\ 1194]  \\
Cost2     & $\mathcal{U}niform$ & [49.6 ,\  74.4] \\
Cost3     & $\mathcal{U}niform$ & [5.6 ,\  8.4 ] \\
Cost4     & $\mathcal{U}niform$ & [4 ,\    8  ]  \\
Duration1 & $\mathcal{U}niform$ & [0.8 ,\  1.2 ] \\
Duration2 & $\mathcal{U}niform$ & [0.8 ,\  1.2]  \\
\hline
\end{tabular}
\caption{VME: distribution of the eighteen parameters.}
\label{tab:VME:2}
\end{table}

In this case, we are particulary interested by the VAN indicator. It is defined by
\begin{eqnarray*}
VAN &=& \PP( G( \xx, \mathbf{W}) \leq 0),
\end{eqnarray*}
where $\mathbf{W} = (W_1,\ldots, W_4)$ and $\xx$ the eighteen fixed parameters.

\begin{figure}
\center
		\includegraphics[width= \textwidth, bb=0 0 1076 587]{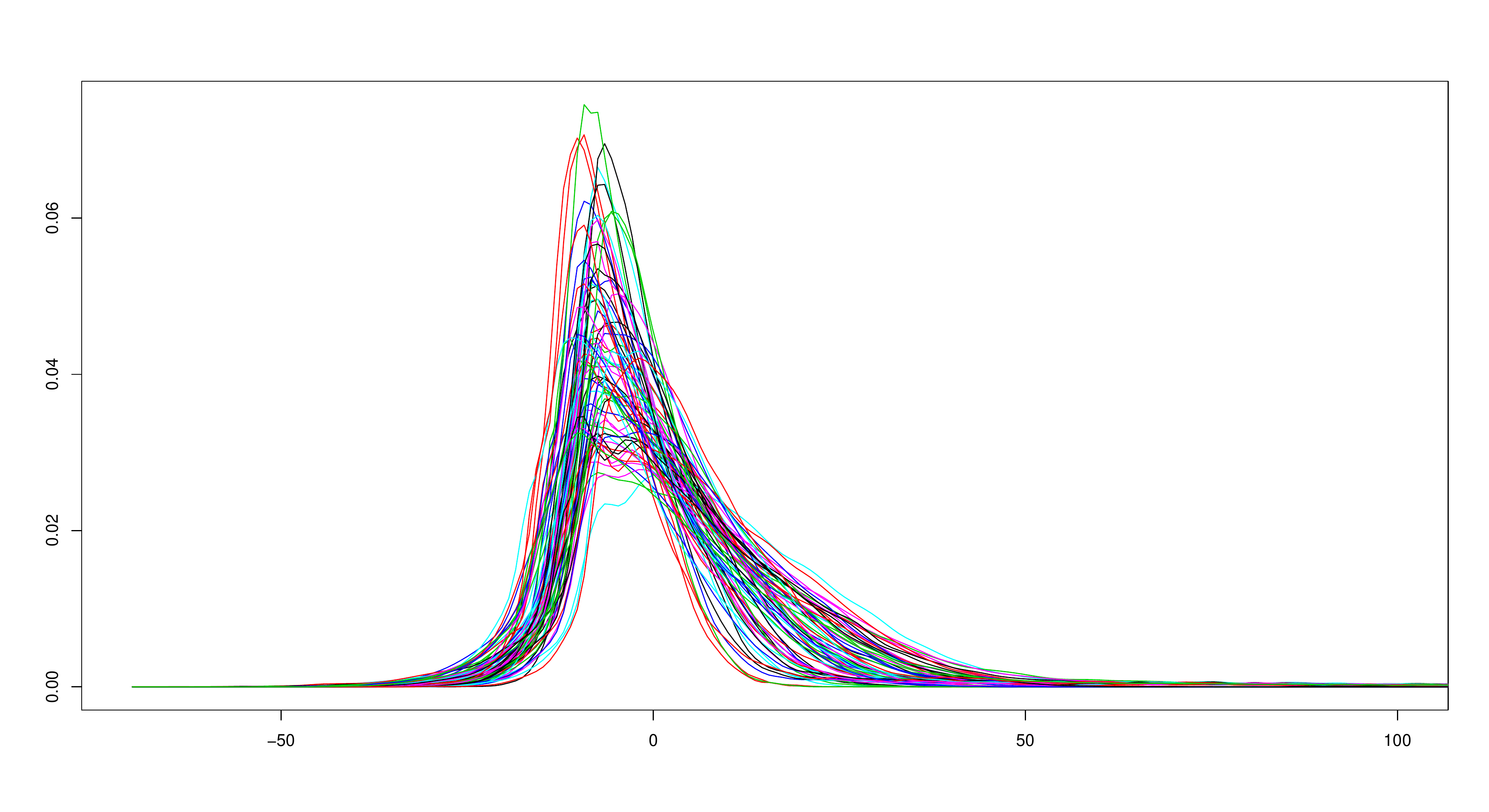}
		\caption{VME: representations of outputs for $N = 75$ different parameters.}
	\label{figure:VME:1:data}
\end{figure}

Figure \ref{figure:VME:1:data} represents 75 outputs from the VME model. The kernel regression is applied to designs of experiments of increasing sizes such that smaller designs are included in greater ones. A test sample of 500 points is used to assess the efficiency of the regression. 

\subsubsection{Kernel regression}

\begin{figure}
\center
	\includegraphics[scale=0.45, bb=0 0 576 576]{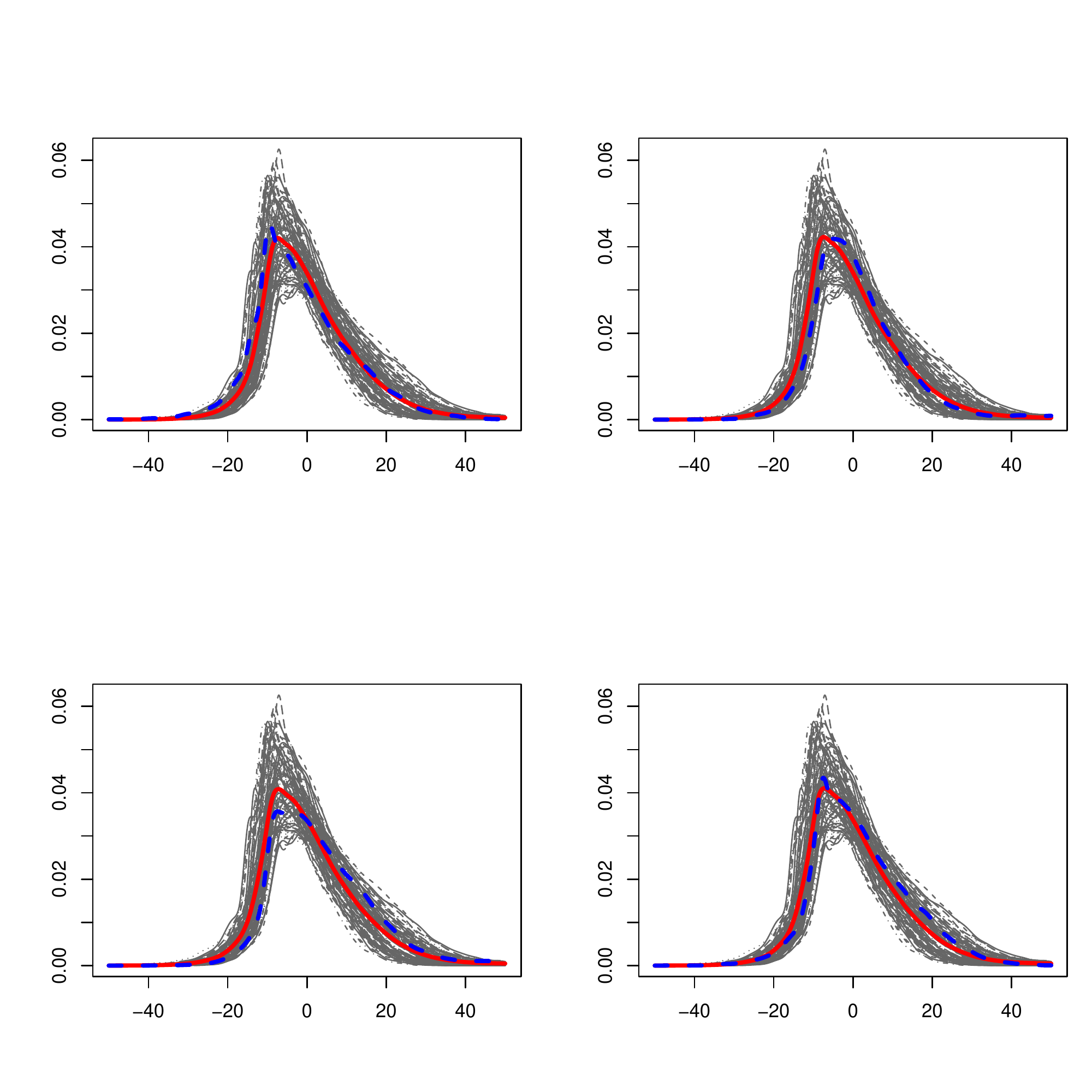}
	\caption{VME: curves in dashed gray lines are probability density functions, red plain line represents the estimation of a density for a fixed $x_0$ and dashed blue line is the true probability density function for $N = 100$.}
\label{figure:VME:4}
\end{figure}

\begin{figure}
\center
		\includegraphics[scale = 0.4, bb=0 0 720 720]{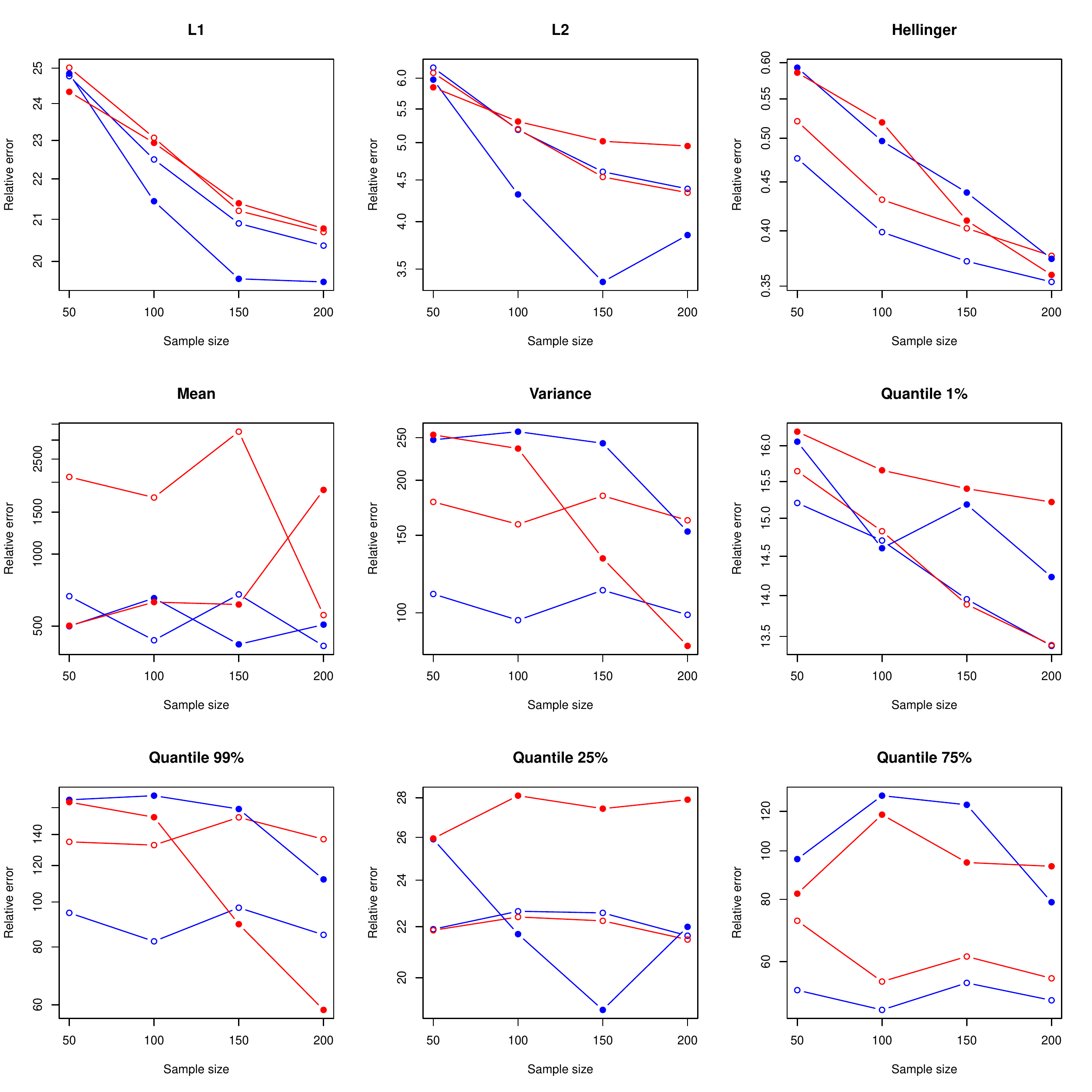}
	\caption{VME: comparison of the relative error for the two estimators with an isotropic and anisotropic (plain cicrle) bandwidth. Blue: estimators given by the $L^2$ norms. Red : estimator given by the Hellinger distance}
\label{figure:VME:err:uni_multi}
\end{figure}

Figure \ref{figure:VME:4} represents for different values of $x_0$ the true probability density function (blue dashed line) and the estimation obtained from the two estimators of the kernel regression (red plain line and orange dots line). In this section, the  kernel regression is tested on this real case study with isotropic and anisotropic bandwidth (Figure \ref{figure:VME:err:uni_multi}). 


In red is represented the relative error for the estimator based on the Hellinger distance, in blue for the estimator based on the $\xLn{2}$ norm. The estimation seems equivalent for the two estimators. The results are equivalent for anisotropic and isotropic bandwidth. For the anisotropic bandwidth, the estimation is less stable and the estimation of $H$ becomes much more time-consuming.

The relative error in norm is smaller for VME (with isotropic and anisotropic bandwidth) than for the toy example 2, whereas the dimension is higher. That shows that the quality of estimation depends not only on the dimension but also on the regularity of the model, if the bandwith is sufficiently well estimated. 

Bad results are found for the variance, 99\% and 99.9\% quantiles. The relative error for the mean is also very high. The real mean is quite close to zero, so that in the relative error computation, the numerator is divided by a quantity close to zero. 

%

\begin{figure}
\center
		\includegraphics[scale = 0.5, bb=0 0 494 356]{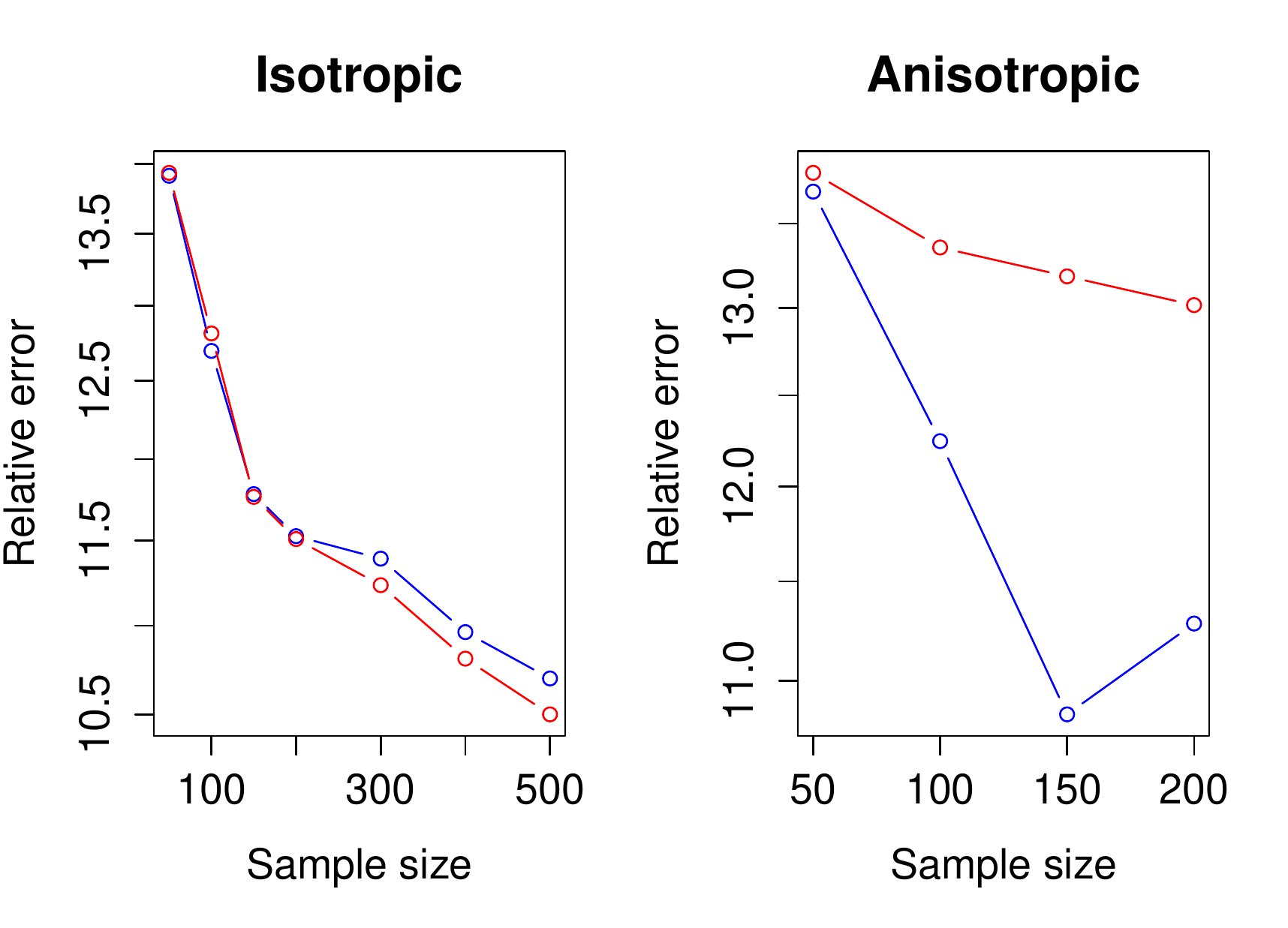}
	\caption{VME: comparison of the relative error for the $\PP(VAN<0)$. Left: isotropic bandwidth. Right: anisotropic bandwidth. Blue: estimators given by the $L^2$ norms. Red : estimator given by the Hellinger distance}
\label{figure:VAN:kernel}
\end{figure}

Figure \ref{figure:VAN:kernel} represents the estimation of the $\PP(VAN<0)$ for isotropic and anisotropic bandwidths. In this case, choosing an anisotropic bandwidth does not bring information.

\subsubsection{Functional decomposition methods}

\begin{figure}
\centering
		\includegraphics[scale = 0.35, bb=0 0 1296 864]{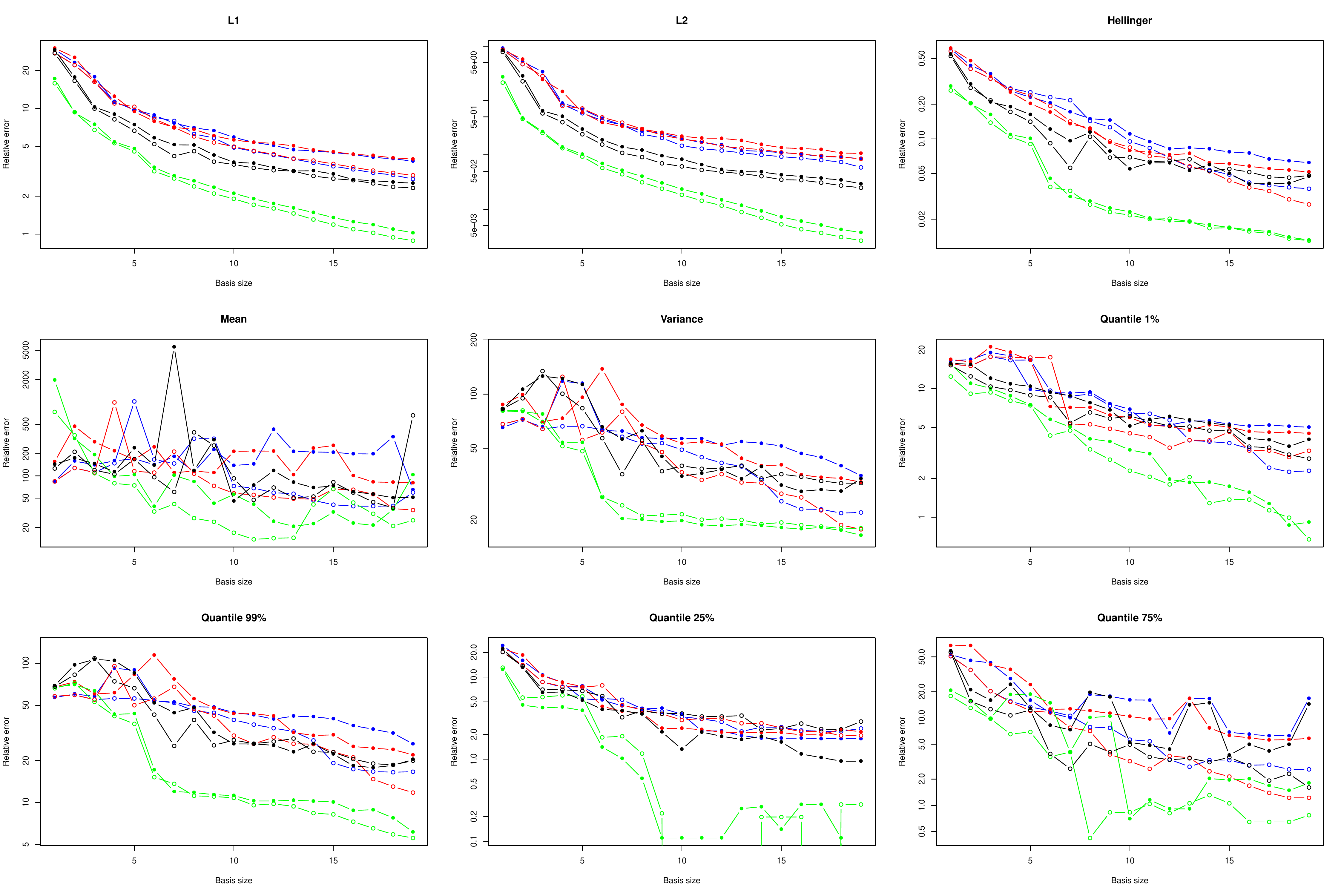}
	\caption{VME: comparison of the relative error for the 9 quantities in function of the size of the basis $q$. Blue: MMP decomposition given by $\xLn{2}$ norm with a learning set size of $N = 50$ (circles) and 100 (filled circles). Red: MMP decomposition given by the Hellinger distance. Black: AQM decomposition. Green: CPCA decomposition.}
\label{fig:vme:decompositions}
\end{figure}

In this section, the four functional decompositions MMP with $\xLn{2}$ and Hellinger distances, AQM and CPCA are compared. Figure \ref{fig:vme:decompositions} represents, in logarithmic scale, the relative errors on the different norms, modes and quantiles versus the basis size (from 1 to 20).
First, the relative errors are low for all quantities of interest, except for the mean and variance. For the variance, the errors are over 20\% with a decomposition basis with 20 functions. The very high errors on the mean can be explained by the fact that the mean to be estimated is close to zero for most of densities in the VME case. The errors on mean, variance have moreover a very chaotic behavior. The errors on distances decrease very quickly. 
For modes and quantiles, the errors of the four methods do not decrease steadily for all quantities of interest. 
CPCA outperforms other methods for the $\xLn{1}$, $\xLn{2}$, Hellinger distances, the variance and the 99\% and 25\% quantiles. The performances of AQM and MMP are very close for modes and quantiles. CPCA seems more stable than MMP and AQM on most of the quantities of interest.

\begin{figure}
\centering
		\includegraphics[scale = 0.35, bb=0 0 720 432]{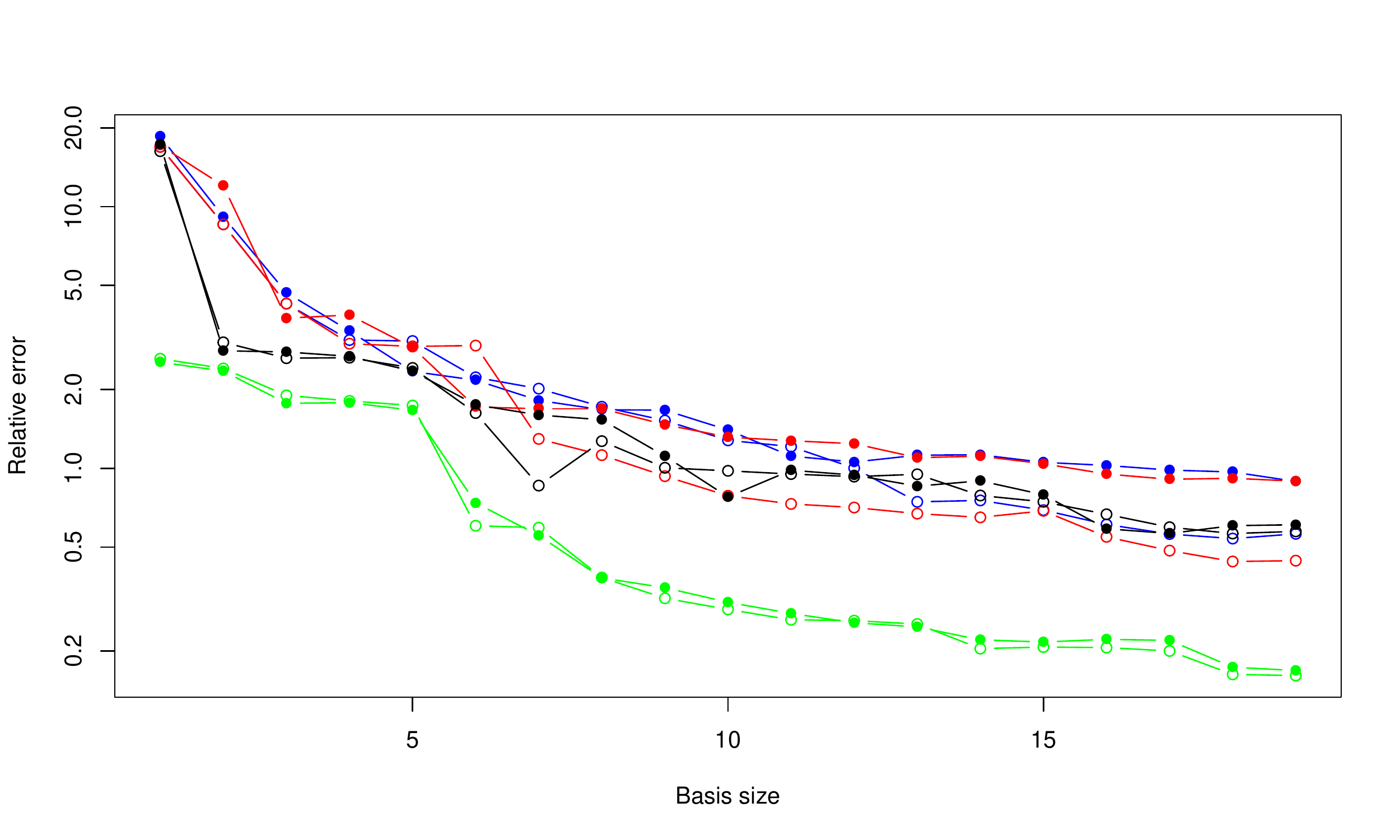}
	\caption{VME: comparison of the relative error for $P(VAN<0)$  in function of the size of the basis $q$. Blue: MMP decomposition given by $\xLn{2}$ norm with a learning set size of $N = 50$ (circles) and 100 (filled circles). Red: MMP decomposition given by the Hellinger distance. Black: AQM decomposition. Green: CPCA decomposition.}
\label{fig:vme:decompositions:van}
\end{figure}

In Figure \ref{fig:vme:decompositions:van}, the errors on the probability for the VAN to be negative $P(VAN<0)$ for the different decomposition methods are represented in function of the basis size. The errors are quite low and decrease quickly. The CPCA method outperforms the others for both learning sets sizes. The other three methods give similar results.

	\subsection{IFPEN application: GIBBS test case}
		\label{section:ifp}
		\subsubsection{Molecular Modelling}

\begin{figure}
	\center
	\includegraphics[scale = 0.1, bb=0 0 1557 603]{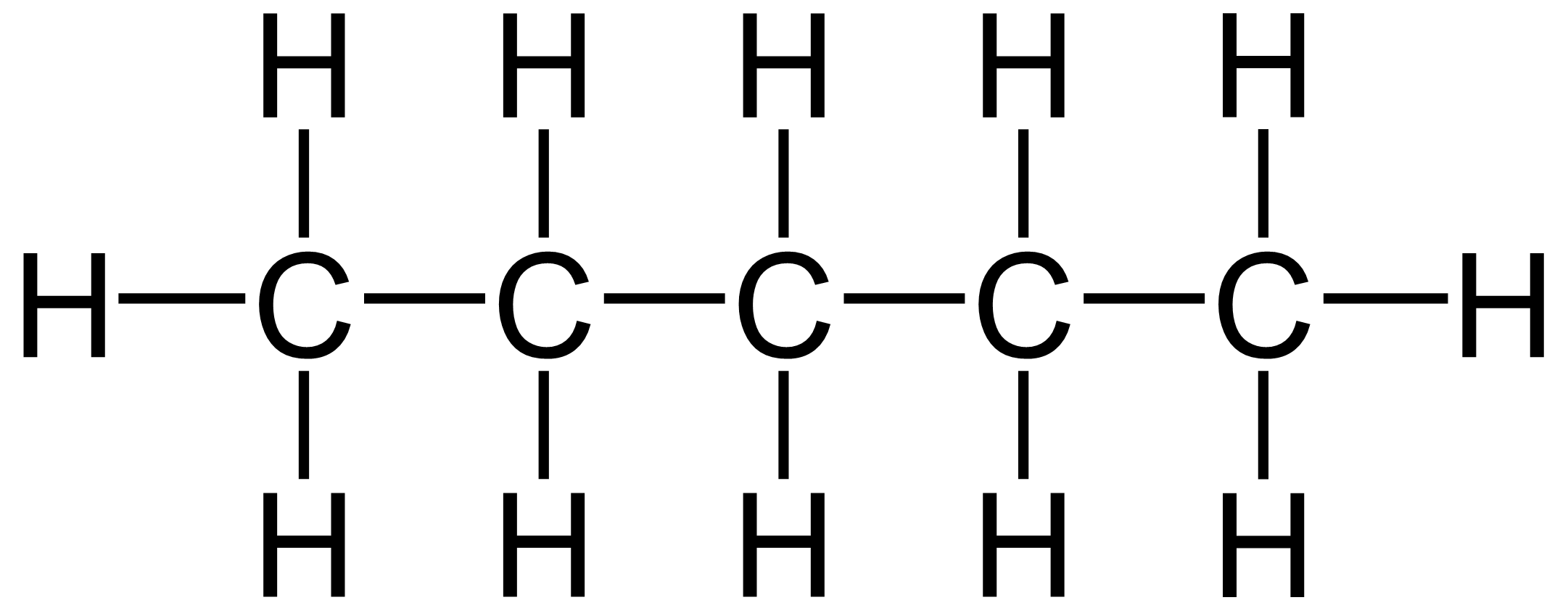}
	\caption{The structural formula of pentane. (Wikipedia)}
	\label{figure:GIBBS:pentane}
\end{figure}

\begin{figure}
\center
\includegraphics[scale=0.4, bb=0 0 1080 432]{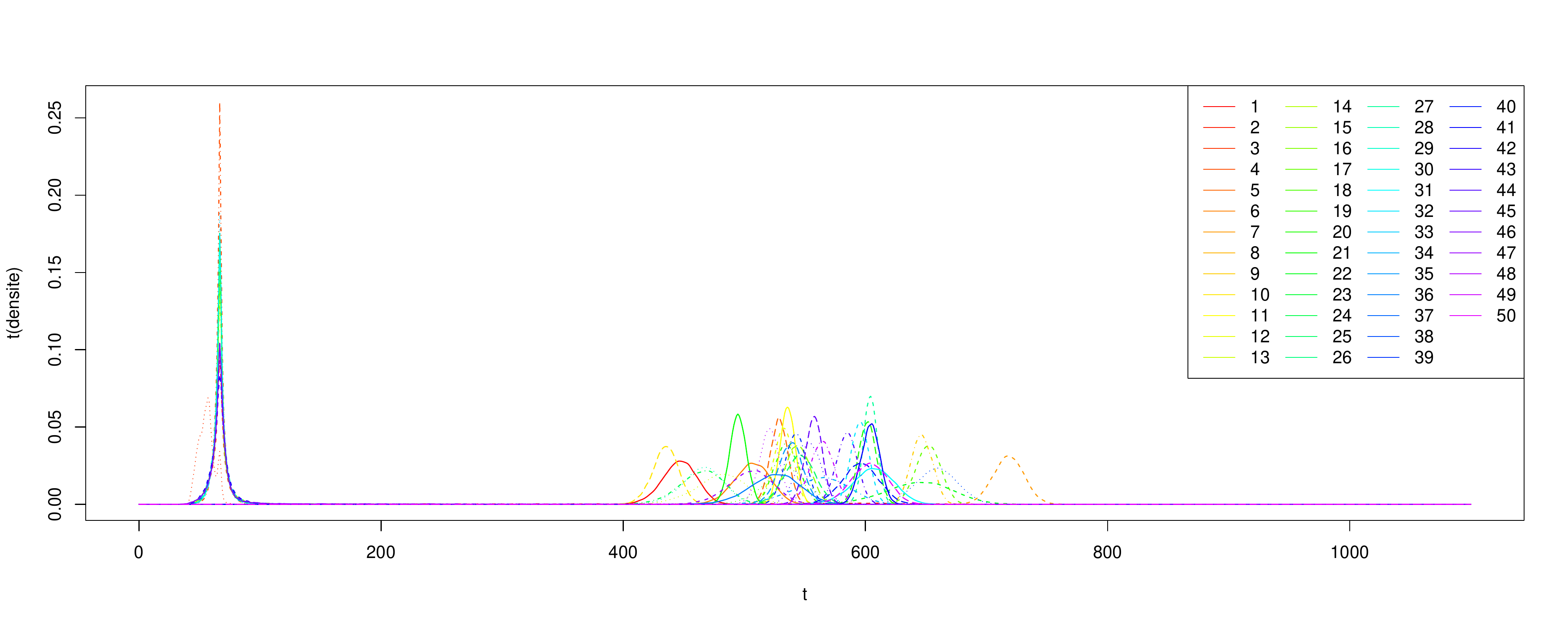}
\caption{GIBBS test case: distributions of the volumetric mass density for different sets of input parameters.}
\label{figure:GIBBS:density:pdf}
\end{figure}

\begin{figure}
\center
\includegraphics[scale=0.4, bb=0 0 1080 432]{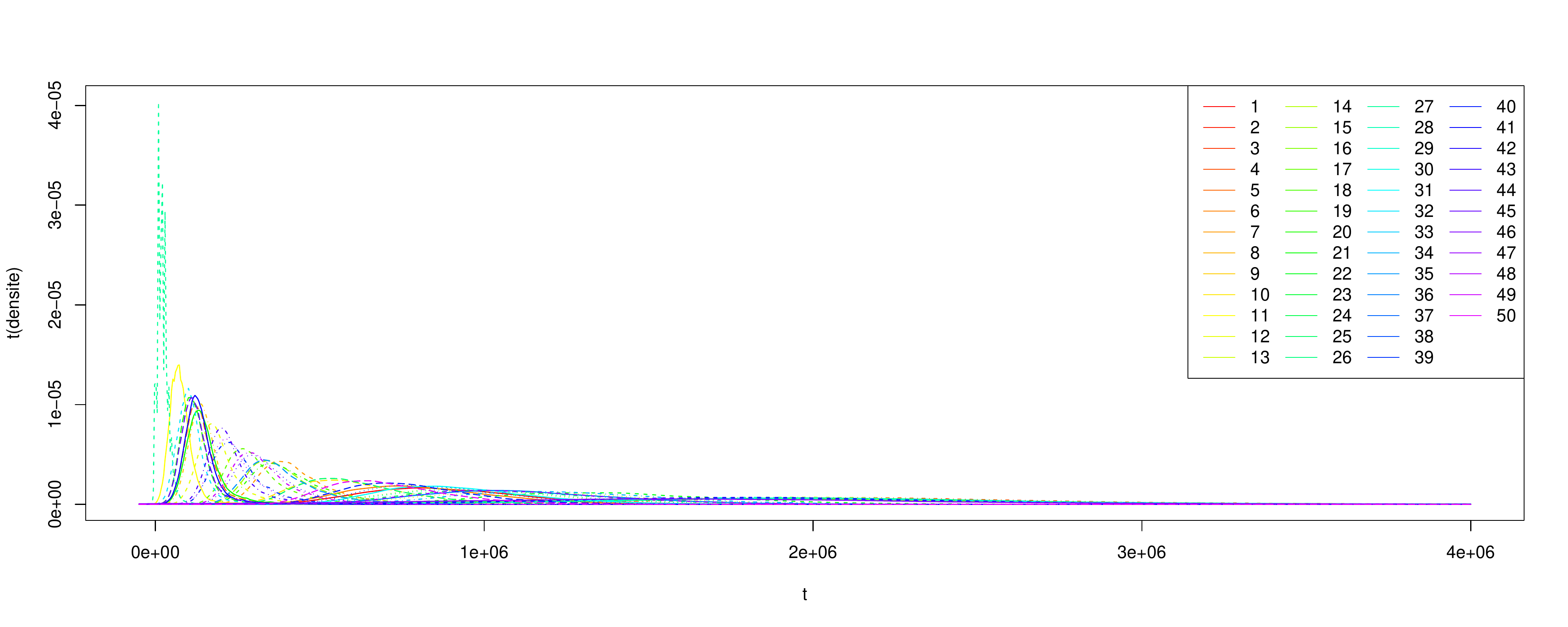}
\caption{GIBBS test case: distributions of the vapor pressure for different sets of input parameters.}
\label{figure:GIBBS:pressure:pdf}
\end{figure}

At IFP Energies Nouvelles, the researchers of the Department of Thermodynamics and Molecular Modelling elaborate models for the structure of molecules such as hydrocarbons and alcohols.
In the following numerical study we consider a system of \emph{pentane} molecules at a \emph{reduced temperature} of $0.75$.
Pentane is an \emph{alkane} with five carbon atoms ($\mathrm{C}_5\mathrm{H}_{12}$, \textit{cf.} Figure \ref{figure:GIBBS:pentane}.)
Here we use the Anisotropic United Atoms model (AUA) \cite{Ungerer00} to represent the molecular structure of pentane.
The carbon ($\mathrm{C}$) and hydrogen ($\mathrm{H}$) atoms in each of the two terminal \emph{methyls} ($\mathrm{CH}_3$) and each of the three \emph{methylene} bridges ($\mathrm{CH}_2$) are treated as a single \emph{interaction center}.
To each of these two types of interaction centers correspond three parameters of the model:
an energy parameter ($\varepsilon_{\mathrm{CH}_3}$ and $\varepsilon_{\mathrm{CH}_2}$), a size parameter ($\sigma_{\mathrm{CH}_3}$ and $\sigma_{\mathrm{CH}_2}$), and a displacement parameter ($\delta_{\mathrm{CH}_3}$ and $\delta_{\mathrm{CH}_2}$). These six parameters will be the components of the set of input parameters of the numerical codes under consideration:
\begin{equation*}
	\xx=\left(\varepsilon_{\mathrm{CH}_3},\varepsilon_{\mathrm{CH}_2},\sigma_{\mathrm{CH}_3},\sigma_{\mathrm{CH}_2},\delta_{\mathrm{CH}_3},\delta_{\mathrm{CH}_2}\right)\in\xR^6.
\end{equation*}\indent
Here the matter of interest is the prediction of two macroscopic properties of a pentane system in \emph{chemical equilibrium}: the \emph{volumetric mass density} $\rho_\mathrm{liq}^\mathrm{eq}$ of the liquid phase and the \emph{vapor pressure} $P_\mathrm{gas}^\mathrm{eq}$ (pressure of the gas phase).
The algorithm carried out to numerically predict $\rho_\mathrm{liq}^\mathrm{eq}$ and $P_\mathrm{gas}^\mathrm{eq}$ -- for a given set of parameters $\xx$ -- relies on the \emph{fundamental postulate of statistical mechanics}: an isolated system in equilibrium is found with equal probability in each of its accessible \emph{microstates}. 
Thus, $\rho_\mathrm{liq}^\mathrm{eq}$ (respectively $P_\mathrm{gas}^\mathrm{eq}$) can be approximated by the average of the values of the volumetric mass density of the liquid phase $\rho_\mathrm{liq}$ (respectively the vapor pressure $P_\mathrm{gas}$) of a large number of accessible microstates of the system.
Therefore the algorithm is essentially a loop, each step consisting of the following two instructions:
\begin{enumerate}[(1)]
	\item randomly generate an accessible microstate of the system (characterized by the positions of the molecules and the volume of the system),
	\item compute the value of $\rho_\mathrm{liq}$ (respectively $P_\mathrm{gas}$) given by the model for this particular microstate.
\end{enumerate}
Let us denote $G_{\rho_\mathrm{liq}}(\xx)$ (respectively $G_{P_\mathrm{gas}}(\xx)$) the probability density function of the random value $\rho_\mathrm{liq}$ (respectively $P_\mathrm{gas}$) computed at each loop step.
The two instructions above are iterated a hundred million times. 
After the execution of the loop, the arithmetical mean of these millions of values is computed -- it is a Monte Carlo approximation of the integral of $G_{\rho_\mathrm{liq}}(\xx)$ (respectively $G_{P_\mathrm{gas}}(\xx)$) -- and is returned by the algorithm as the prediction of $\rho_\mathrm{liq}$ (respectively $P_\mathrm{gas}$.) 
The whole algorithm is carried out in about twenty-four hours on a supercomputer (for each set $\xx$ of input parameters.)\newline\indent
In order to circumvent this time-consuming process, we want to build a metamodel of the numerical codes mapping the input parameters $\xx$ to the probability density functions $G_{P_\mathrm{gas}}(\xx)$ and $G_{\rho_\mathrm{liq}}(\xx)$ of the random outputs. We had at our disposal a sample of size $N=50$. Figure \ref{figure:GIBBS:density:pdf} and Figure \ref{figure:GIBBS:pressure:pdf} show the corresponding $50$ probability functions for volumetric mass density and vapor pressure respectively.

%

\subsubsection{Kernel regression}

We applied the kernel regression method on the numerical codes providing the volumetric mass density of the liquid phase and the vapor pressure. 
The results provided by the $\xLn{2}$ estimator with isotropic bandwidth are provided in Table \ref{table:GIBBS:KR:massvol:err} for the volumetric mass density and in Table \ref{table:GIBBS:KR:pressure:err} for the vapor pressure.
The errors have been computed by the Leave-One-Out methodology as in Section \ref{castem:kr}.
We can see that most of the results are disappointing: all relative errors are higher than $30\%$, excepted -- surprisingly -- the relative Hellinger error of the approximation. 
The latter amounts to nearly $7\%$ for the volumetric mass density, and only nearly $0.05\%$ for the vapor pressure.
The errors for quantile 1\% and quantile 25\% are particularly high.
The Hellinger estimator give almost the same results: the difference between the error for the $\xLn{2}$ estimator and the error for the Hellinger estimator is lower than $10^{-5}$ in every case. 
That is why they are not charted here.

\begin{table}
\center
\begin{tabular}{ccc}
\hline
L1  & L2 & Hellinger\\
151.1\% & 375.9\% & 7.089\% \\
\hline \\
\hline \\
Mean & Variance & Quantile 1\%\\
50.77\% & 98.07\% & 513.8\% \\
\hline \\
\hline \\
Quantile 99\% & Quantile 25\% & Quantile 75\%\\
30.44\% & 198.4\% & 41.52\% \\
\hline 
\end{tabular}
\caption{GIBBS test case: kernel regression relative errors of the $\xLn{2}$ estimator for the volumetric mass density.}
\label{table:GIBBS:KR:massvol:err}
\end{table}

\begin{table}
\center
\begin{tabular}{ccc}
\hline
L1  & L2 & Hellinger\\
123.7\% & 254.1\% & 0.04878\% \\
\hline \\
\hline \\
Mean & Variance & Quantile 1\%\\
69.86\% & 80.72\% & 667.63\% \\
\hline \\
\hline \\
Quantile 99\% & Quantile 25\% & Quantile 75\%\\
74.82\% & 156.1\% & 73.5\%\\
\hline 
\end{tabular}
\caption{GIBBS test case: KR relative errors of the $\xLn{2}$ estimator for the vapor pressure.}
\label{table:GIBBS:KR:pressure:err}
\end{table}

\begin{figure}
\center
\includegraphics[scale = 0.4, bb=0 0 720 576]{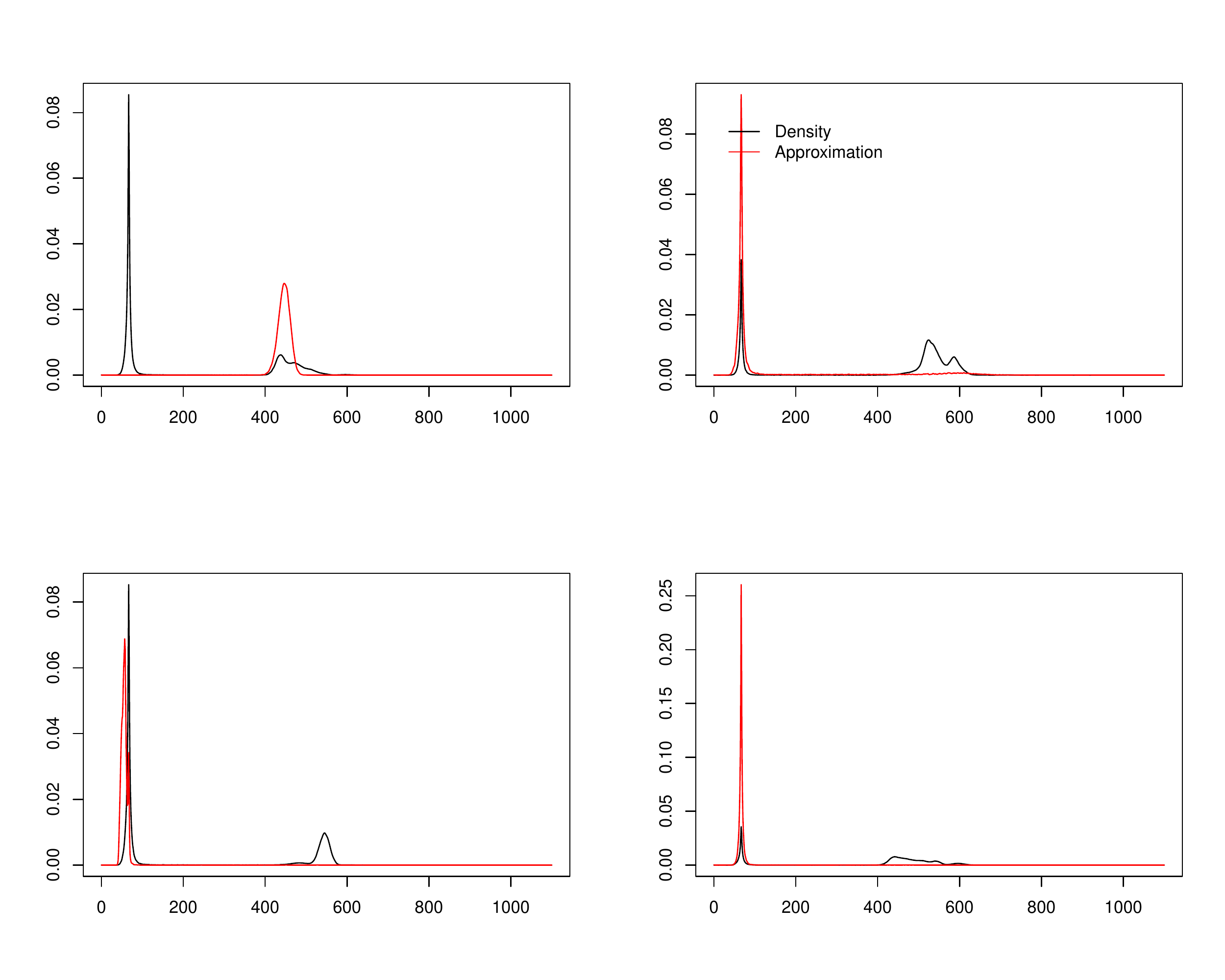}
\caption{GIBBS test case: distribution of the volumetric mass density and its KR approximation for the sets of input parameters 1, 2, 3 and 4.}
\label{figure:GIBBS:KR:massvol:plot}
\includegraphics[scale = 0.4, bb=0 0 720 576]{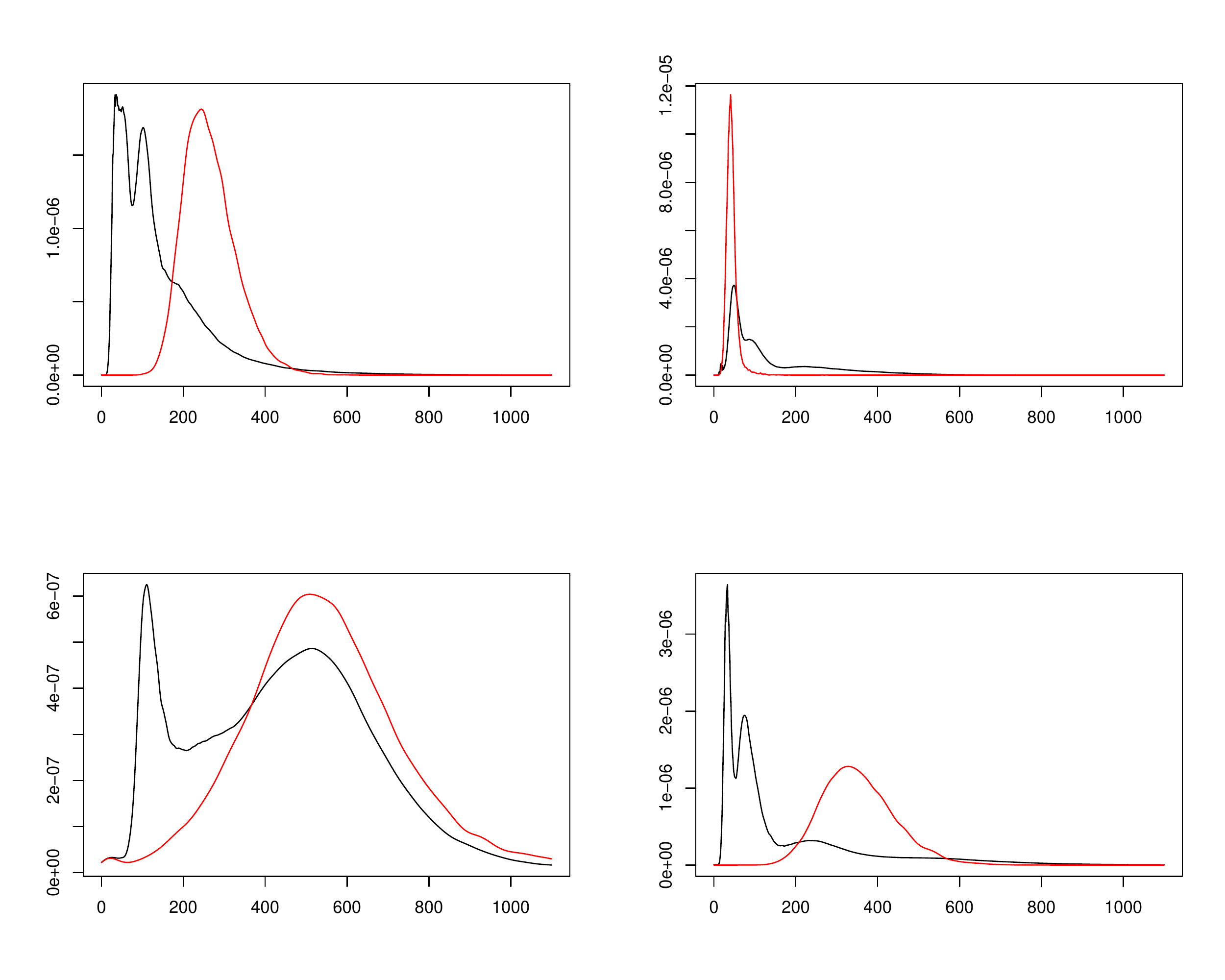}
\caption{GIBBS test case: distribution of the vapor pressure and its KR approximation for the sets of input parameters 31, 32, 33 and 34.}
\label{figure:GIBBS:KR:pressure:plot}
\end{figure}

In Figure \ref{figure:GIBBS:KR:massvol:plot} and Figure \ref{figure:GIBBS:KR:pressure:plot} we plot the probability density functions of the volumetric mass density and the vapor pressure with their kernel regression approximations, for arbitrarily chosen sets of input parameters. 
We can see how poorly the densities are approximated.

\subsubsection{Functional decomposition methods}

\begin{figure}
\center
		\includegraphics[scale=0.35, bb=0 0 1296 864]{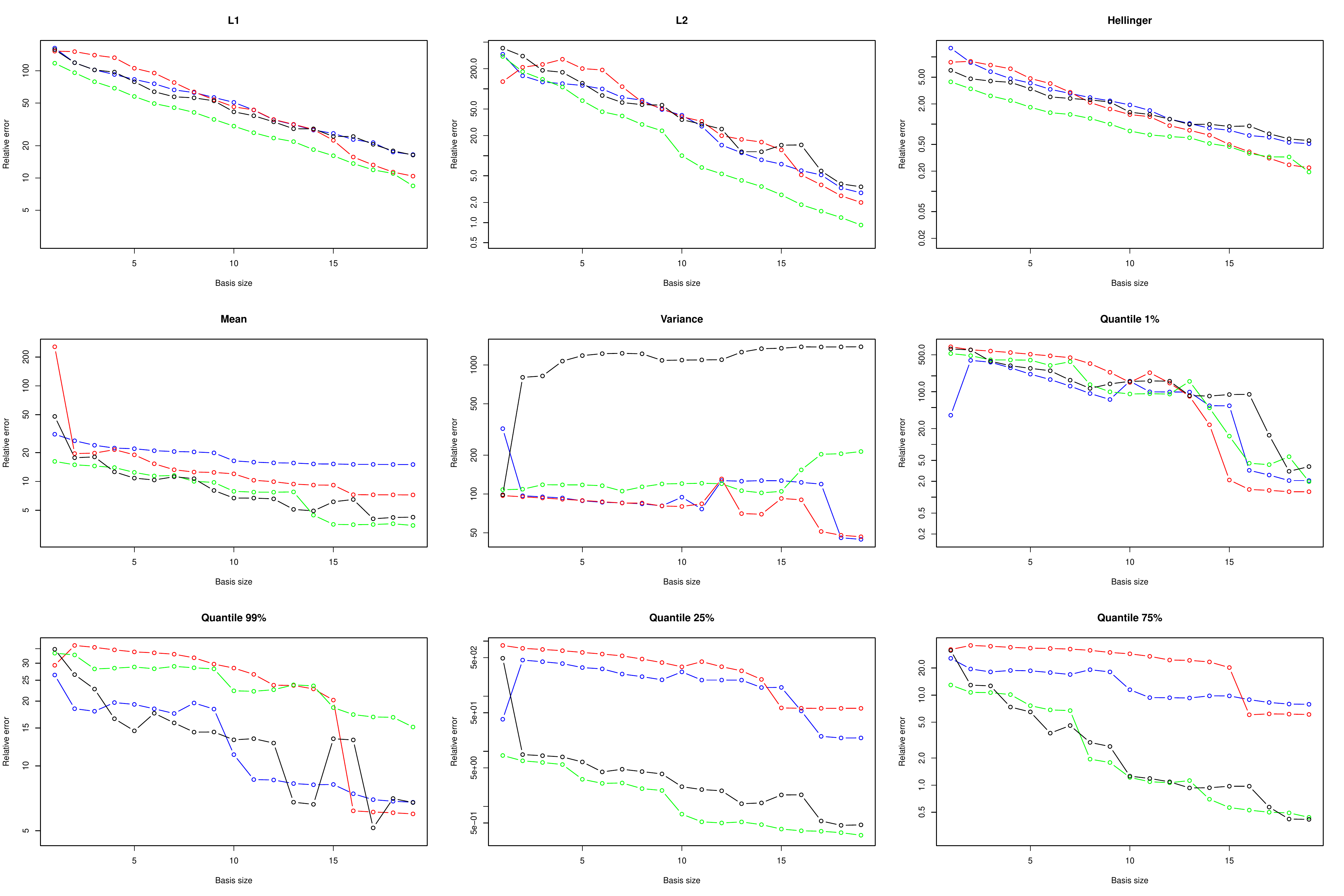}
	\caption{GIBBS test case: comparison of the relative errors on different quantities for the volumetric mass density in function of the size of the basis $q$. Blue: MMP decomposition given by $\xLn{2}$ norm. Red: MMP decomposition given by the Hellinger distance. Black: AQM decomposition. Green: CPCA decomposition.}
\label{figure:massvol:decompositions}
\end{figure}

In this section, the four functional decompositions MMP with $\xLn{2}$ and Hellinger distances, AQM and CPCA are compared. Figure \ref{figure:massvol:decompositions} represents, in logarithmic scale, the relative errors on the different norms, modes and quantiles versus the basis size (from 1 to 20) for the volumetric mass density of the liquid phase.
The error curves are globally decreasing excepted the one corresponding to the variance relative error for AQM. 
Moreover, the relative error on variance is high even for basis with 20 functions and especially for AQM. 
The behavior of the AQM approximation is more erratic for the 99\% quantile.
Both MMP methods perform worse on the 25\% and 75\% quantiles than other methods.
Other relative errors are low. 
CPCA gives better results for $\xLn{1}$, $\xLn{2}$, Hellinger distances and the 25\% quantile. For 25\% and 75\% quantiles, errors on both MMP bases are higher than for others methods.

\begin{figure}
\center
		\includegraphics[scale=0.35, bb=0 0 1296 864]{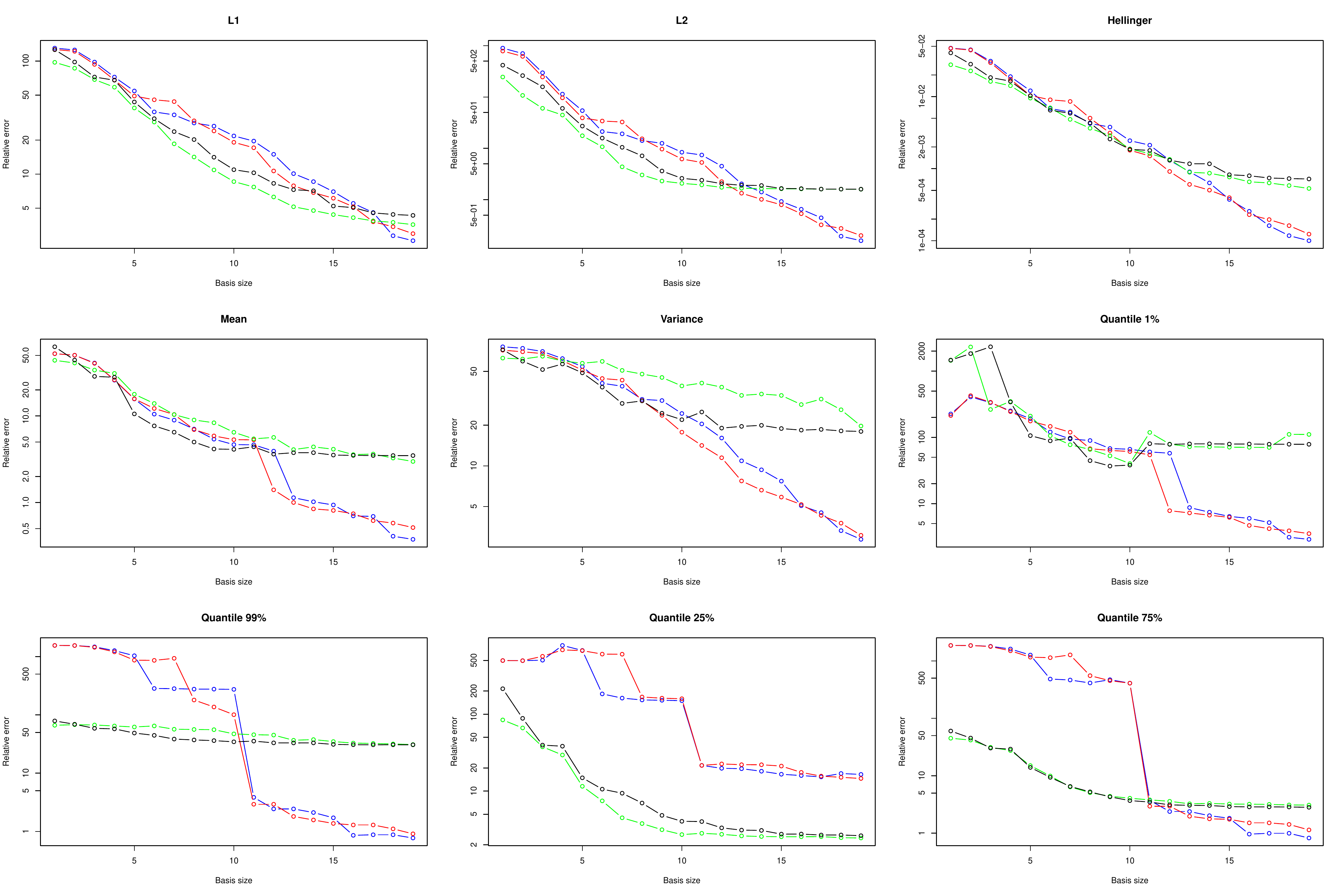}
	\caption{GIBBS test case: comparison of the relative errors on different quantities for the vapor pressure in function of the size of the basis $q$. Blue: MMP decomposition given by $\xLn{2}$ norm. Red: MMP decomposition given by the Hellinger distance. Black: AQM decomposition. Green: CPCA decomposition.}
\label{figure:pressure:decompositions}
\end{figure}

The results for the vapor pressure are presented in Figure \ref{figure:pressure:decompositions}. The errors on all quantities of interest are more steady than for the volumetric mass. In this case, CPCA and AQM have very different behavior compared to both MMP methods. For all quantities of interest except $\xLn{1}$, CPCA and AQM errors have a plateau for higher basis sizes. For higher basis sizes, MMP method performs better than AQM and CPCA for all studied quantities, except the 25\% quantile. Thus, contrary to what happens in the other presented numerical examples, CPCA method does not perform better than AQM and MMP for distances. 
For MMP decompositions, the errors on the mean and quantiles display a jump for basis sizes of 12 or 13.

\section{Conclusion}
	The aim of this work was to build a metamodel for code outputs which are probability density functions. A first idea to design such a model was to build a convex combination of the sample probability density functions. To this mean, we first adapted the well-known kernel regression developed by \cite{nadaraya1964} and \cite{watson1964}. 
Second, we proposed to approximate the sample probability density functions on a functional basis in order to reduce the problem dimension. In this way, the probability density functions are characterized by their coefficients on the basis, so that the problem becomes finite-dimensional. 
Two methods have been proposed to build a functional basis. The first method is adapted from Magic Points Method \cite{Maday2007} and builds iteratively the basis by adding the sample function that maximizes the approximation error on the previous basis. The second one, called Alternate Quadratic Minimization (AQM), aims at minimizing the $\xLn{2}$ approximation error under the defined constraints. Both methods require the resolution of constrained optimization problems.
Then, in future work, a metamodel could be adjusted on the coefficients, to link them with the uncertain inputs, and provide a global metamodel for probability density function outputs.

Our methods have been tested on two analytical test cases and three industrial applications proposed by CEA, EDF and IFPEN. The kernel regression-based metamodel performs well on the analytical test cases, but shows its limits on the three numerical codes. Indeed, kernel regression is efficient in low dimension, but the difficulty to estimate the parameters increases with the dimension and the size of the sample. Moreover, this method is not adapted for numerical codes for which the influence of the controllable parameters is much more important than one of the uncontrollable parameters.
Both Modified Magic Points (MMP) and AQM methods give good results on all the test cases. However, in most cases, the method proposed by Kneip and Utikal \cite{Kneip2001} performs better than the proposed ones.


Metamodels linking the uncertain inputs of the code and the functional basis coefficients remains to be defined. The main difficulty to build this metamodel lies in ensuring that the sum of the functional basis coefficients is equal to one.
One of the limits of our methods is that all quantities of interest are not as well approximated. In particular, quantiles are often poorly approximated. Indeed, all the proposed methods are based on norm minimization, and the norms used are not well suited for the study of extreme values. It could be interesting to adapt the proposed methods to norms giving more weights on the tails of the densities. The use of the Wasserstein metric between probability measures will be also studied.

\begin{acknowledgement}
We would like to thank our supervisors Nicolas Bousquet (EDF), Fr\'ed\'eric Delbos (IFPEN), Bertrand Iooss (EDF), J\'er\^ome Lonchampt (EDF), Rafael Lugo (IFPEN), Amandine Marrel (CEA) and Nadia P\'erot (CEA).
We would also like to thank Anthony Nouy (\'Ecole Centrale de Nantes) and S\'ebastien Da Veiga (IFPEN) for their guidance.
Furthermore, we would like to thank Nicolas Champagnat (Institut \'Elie Cartan de Lorraine) and Tony Leli\`evre (\'Ecole des Ponts ParisTech) who organized the CEMRACS 2013 in collaboration with Anthony Nouy.
Finally, we would like to thank CEA, EDF and IFPEN for funding this research project.
\end{acknowledgement}

\bibliographystyle{plain}
\bibliography{cemracs}

\appendix
\section{Appendix}
	\subsection{$\xLn{2}$ distance-based kernel estimator}
	\label{section:L2.kernel.equivalence}
	We prove here the equivalence between the classical kernel estimator and the $\xLn{2}$ distance-based kernel estimator introduced in Sections \ref{section:KernelRegression} and \ref{section:KernelRegressionHellinger} respectively.
	Indeed the classical kernel estimator \eqref{eq:KR} can be retrieved by writing the optimality condition of the minimization problem \eqref{eq:KR2} defining the $\xLn{2}$ distance-based kernel estimator:
	\begin{align*}
		& \quad\sum_{i=1}^N-2K_H(\xx_i,\xx_0)\int\limits_I\left(f_i-\hat{f}_0\right)=0
		\\\Longleftrightarrow & \quad\int\limits_I\left(\sum_{i=1}^NK_H(\xx_i,\xx_0)\left(f_i-\hat{f}_0\right)\right)=0
		\\\Longleftrightarrow & \quad\sum_{i=1}^NK_H(\xx_i,\xx_0)\left(f_i-\hat{f}_0\right)=0
		\\\Longleftrightarrow & \quad\hat{f}_0=\sum_{i=1}^N\Frac{K_H(\xx_i,\xx_0)}{\sum\limits_{j=1}^N K_H(\xx_j,\xx_0)}f_i.
	\end{align*}

\subsection{Hellinger distance-based kernel estimator}
	\label{section:Hellinger.kernel.equivalence}
	We prove here the equivalence between expressions \eqref{eq:Hellinger:2} and \eqref{eq:Hellinger:1} of the Hellinger distance-based kernel estimator given in Section \ref{section:KernelRegressionHellinger}.
	The constraint that the integral of $\hat{f}_0$ is 1 in the optimization problem \eqref{eq:Hellinger:2} is handled using the associated Lagrangian function.
	The problem becomes
	\begin{equation*}
		\hat{f}_0=\argmin_{f,\lambda}L(f,\lambda),
	\end{equation*}
	with
	\begin{equation*}
		L(f,\lambda)=\sum_{i=1}^N K_H(\xx_i,\xx_0)\int\limits_I\left(\sqrt{f_i}-\sqrt{f}\right)^2-\lambda\left(\int\limits_If-1\right). 
	\end{equation*}
	The first order optimality conditions are
	\begin{equation}
		\label{eq:derivative}
		\left\{\begin{aligned}
			0 & =\Frac{\partial L}{\partial f} =\sum\limits_{i=1}^NK_H(\xx_i,\xx_0)\int\limits_I\left(1-\Frac{\sqrt{f_i}}{\sqrt{\hat{f}_0}}\right)-\lambda\int\limits_I1,
			\\0 & =\Frac{\partial L}{\partial\lambda}=1-\int\limits_I\hat{f}_0.
		\end{aligned}\right.
	\end{equation}
	The first condition in \eqref{eq:derivative} gives
	\begin{align}
		\label{eq:firstcondition}
		& \quad\int\limits_I\left(\sum\limits_{i=1}^NK_H(\xx_i,\xx_0)\left(1-\Frac{\sqrt{f_i}}{\sqrt{\hat{f}_0}}\right)-\lambda\right)=0\nonumber
		\\\Longleftrightarrow & \quad\sum\limits_{i=1}^N K_H(\xx_i,\xx_0)\left(1-\Frac{\sqrt{f_i}}{\sqrt{\hat{f}_0}}\right)-\lambda=0\nonumber
		\\\Longleftrightarrow & \quad\hat{f}_0=\left(\sum\limits_{i=1}^N \Frac{K_H(\xx_i,\xx_0)\sqrt{f_i}}{\sum\limits_{j=1}^N K_H(\xx_j,\xx_0)-\lambda}\right)^2=\Frac{\displaystyle\left(\sum\limits_{i=1}^NK_H(\xx_i,\xx_0)\sqrt{f_i}\right)^2}{\displaystyle\left(\sum\limits_{j=1}^N K_H(\xx_j,\xx_0)-\lambda\right)^2}.
	\end{align}
	Then we replace $\hat{f}_0$ in the second condition of \eqref{eq:derivative} by the expression given in \eqref{eq:firstcondition}:
	\begin{align}
		\label{eq:secondcondition}
\left(\sum\limits_{j=1}^N K_H(\xx_j,\xx_0)-\lambda\right)^2=\int\limits_I\left(\sum\limits_{i=1}^N K_H(\xx_i,\xx_0)\sqrt{f_i}\right)^2.
	\end{align}
	The combined equations \eqref{eq:firstcondition} and \eqref{eq:secondcondition} entail the formula given by \eqref{eq:Hellinger:1} for the kernel estimator associated to the Hellinger distance:
	\begin{equation*}
		\hat{f}_0=\Frac{\displaystyle\left(\sum\limits_{i=1}^NK_H(\xx_i,\xx_0)\sqrt{f_i}\right)^2}{\displaystyle\int\limits_I\left(\sum\limits_{i=1}^N K_H(\xx_i,\xx_0)\sqrt{f_i}\right)^2}.
	\end{equation*}

\end{document}